\documentclass[leqno,12pt]{article}
\usepackage{epsfig,amsfonts,latexsym}

\font\refit=cmcsc10
\font\capit=cmcsc10 
\font\addressit=cmcsc8
\font\eightrm=cmr8

\font\twelvesymb msbm10 at 12pt
 at 14.4truept

\makeatletter                                             
\renewcommand{\@begintheorem}[2]{                        
\sl \trivlist \item                                      
[\hskip \labelsep {\capit #2\ \ #1.}]  }                 
    
\makeatother         
\makeatletter
\def\section{\@startsection{section}{1}{\z@}{-3.5ex plus -1ex minus
 -.2ex}{1.5ex plus .2ex}{\large\bf}}

\def\subsection{\@startsection{subsection}{2}{\z@}{-3.25ex plus -1ex
minus -.2ex}{1.5ex plus .2ex}{\normalsize\it}}

\def\textindent#1{\indent\llap{#1\enspace}\ignorespaces}
\def\itemitem{\par\indent \hangindent2\parindent \textindent}

\newcommand{\numberequationsassubsubsections}

\newtheorem{prop}{Proposition}[subsection]
\newtheorem{lemm}[prop]{Lemma}
\newtheorem{theo}[prop]{Theorem}

\newtheorem{atiyah-bott}[prop]{Atiyah-Bott fixed point theorem}
\newtheorem{define}[prop]{Definition}
\newtheorem{lemm-def}[prop]{Lemma-Definition}

\newtheorem{rem}[prop]{\it Remark}
\newtheorem{rems}[prop]{\it Remarks}
\newtheorem{notat}[prop]{\it Notation}

\def\statement #1. #2\par{\medbreak\noindent{\bf#1.\enspace}{\sl#2}\par
	\ifdim\lastskip<\medskipamount \removelastskip\penalty55\medskip\fi}

\font\dynkfont=cmsy10 scaled\magstep4    \skewchar\dynkfont='60
\def\dynk{\textfont2=\dynkfont}
\def\hr#1,#2;{\dimen0=.4pt\advance\dimen0by-#2pt
              \vrule width#1pt height#2pt depth\dimen0}
\def\vr#1,#2;{\vrule height#1pt depth#2pt}
\def\blb#1#2#3#4#5
            {\hbox{\ifnum#2=0\hskip11.5pt
                   \else\ifnum#2=1\hr13,5;\hskip-1.5pt
                   \else\ifnum#2=2\hr13.5,6.5;\hskip-13.5pt
                                  \hr13.5,3.5;\hskip-2pt
                   \else\ifnum#2=3\hr13.7,8;\hskip-13.7pt
                                  \hr13,5;\hskip-13pt
                                  \hr13.7,2;\hskip-2.2pt\fi\fi\fi\fi
                   $#1$
                   \ifnum#4=0
                   \else\ifnum#4=1\hskip-9.2pt\vr22,-9;\hskip8.8pt
                   \else\ifnum#4=2\hskip-10.9pt\vr22,-8.75;\hskip3pt
                                  \vr22,-8.75;\hskip7.1pt
                   \else\ifnum#4=3\hskip-12.6pt\vr22,-8.5;\hskip3pt
                                  \vr22,-9;\hskip3pt
                                  \vr22,-8.5;\hskip5.4pt\fi\fi\fi\fi
                   \ifnum#5=0
                   \else\ifnum#5=1\hskip-9.2pt\vr1,12;\hskip8.8pt
                   \else\ifnum#5=2\hskip-10.9pt\vr1.25,12;\hskip3pt
                                  \vr1.25,12;\hskip7.1pt
                   \else\ifnum#5=3\hskip-12.6pt\vr1.5,12;\hskip3pt
                                  \vr1,12;\hskip3pt
                                  \vr1.5,12;\hskip5.4pt\fi\fi\fi\fi
                   \ifnum#3=0\hskip8pt
                   \else\ifnum#3=1\hskip-5pt\hr13,5;
                   \else\ifnum#3=2\hskip-5.5pt\hr13.5,6.5;
                                  \hskip-13.5pt\hr13.5,3.5;
                   \else\ifnum#3=3\hskip-5.7pt\hr13.7,8;
                                  \hskip-13pt\hr13,5;
                                  \hskip-13.7pt\hr13.7,2;\fi\fi\fi\fi
 
                   }}
\def\blob#1#2#3#4#5#6#7{\hbox
{$\displaystyle\mathop{\blb#1#2#3#4#5 }_{#6}\sp{#7}$}}
\def\up#1#2{\dimen1=33pt\multiply\dimen1by#1\hbox{\raise\dimen1\rlap{#2}}}
\def\uph#1#2{\dimen1=17.5pt\multiply\dimen1by#1\hbox{\raise\dimen1\rlap{#2}}}
\def\dn#1#2{\dimen1=33pt\multiply\dimen1by#1\hbox{\lower\dimen1\rlap{#2}}}
\def\dnh#1#2{\dimen1=17.5pt\multiply\dimen1by#1\hbox{\lower\dimen1\rlap{#2}}}
 
\def\rlbl#1{\kern-8pt\raise3pt\hbox{$\scriptstyle #1$}}
\def\llbl#1{\raise3pt\llap{\hbox{$\scriptstyle #1$\kern-8pt}}}
\def\elbl#1{\kern3pt\lower4.5pt\hbox{$\scriptstyle #1$}}
\def\lelbl#1{\rlap{\hbox{\kern-9pt\raise2.5pt\hbox{{$\scriptstyle #1$}}}}}

\def\whtd#1#2#3#4#5{\blob\circ#1#2#3#4{#5}{}}
\def\blkd#1#2#3#4#5{\blob\bullet#1#2#3#4{#5}{}}
\def\whtu#1#2#3#4#5{\blob\circ#1#2#3#4{}{#5}}
\def\blku#1#2#3#4#5{\blob\bullet#1#2#3#4{}{#5}}
\def\whtr#1#2#3#4#5{\blob\circ#1#2#3#4{}{}\rlbl{#5}}
\def\blkr#1#2#3#4#5{\blob\bullet#1#2#3#4{}{}\rlbl{#5}}

\def\rwng{\hbox{$\vbox{\offinterlineskip{
  \hbox{\phantom{}\kern6pt{$\circ$}}\kern-2.5pt\hbox{$\Biggr/$}\kern-0.5pt
  \hbox{\phantom{}\kern-5pt$\circ$}\kern-3.0pt\hbox{$\Biggr\backslash$}
  \kern-1.5pt\hbox{\phantom{}\kern6pt{$\circ$}} }}$}}
 
\def\lwng{\hbox{$\vbox{\offinterlineskip{ \hbox{$\circ$}
  \kern-3.0pt\hbox{\phantom{}\kern6.0pt{$\Biggr\backslash$}}
  \kern-0.5pt\hbox{\phantom{}\kern11pt{$\circ$}}\kern-3.5pt
  \hbox{\phantom{}\kern5.0pt {$\Biggr/$}}\kern-1.0pt\hbox{$\circ$} }}$}}

\def\drwng#1#2#3{\hbox{$\vcenter{ \offinterlineskip{
  \hbox{\phantom{}\kern7pt{$\circ^{\elbl{#3}}$}}
  \kern-2.5pt\hbox{$\Biggr/$}\kern-0.5pt
  \hbox{\phantom{}\kern-5pt$\circ^{ \elbl{#1}}$}
  \kern-3.0pt\hbox{$\Biggr\backslash$}
  \kern-1.5pt\hbox{\phantom{}\kern7pt{$\circ^{\elbl{#2}}$}}  } }$}}
\def\drwngt#1#2#3{\hbox{$\vcenter{ \offinterlineskip{
  \hbox{\phantom{}\kern7pt{$\bullet^{\elbl{#3}}$}}
  \kern-2.5pt\hbox{$\Biggr/$}\kern-0.5pt
  \hbox{\phantom{}\kern-5pt$\circ^{ \elbl{#1}}$}
  \kern-3.0pt\hbox{$\Biggr\backslash$}
  \kern-1.5pt\hbox{\phantom{}\kern7pt{$\bullet^{\elbl{#2}}$}}  } }$}}
 
\def\dlwng#1#2#3{\hbox{$\vcenter{\offinterlineskip{ \hbox{$\lelbl{#1}\circ$}
  \kern-3.0pt\hbox{\phantom{}\kern6.0pt{$\Biggr\backslash$}}
  \kern-0.5pt\hbox{\phantom{}\kern11pt{$\lelbl{#2}\circ$}}\kern-3.5pt
  \hbox{\phantom{}\kern5.0pt {$\Biggr/$}}\kern-1.0pt\hbox{$\lelbl{#3}\circ$}}}$}
 }

\def\rde#1#2#3{\hbox{\phantom{}\kern-4pt\hbox{$\vcenter{\offinterlineskip  \hbox
{
               \raise 4.5pt\hbox{\vrule height0.4pt width13pt depth0pt}
                \kern-1pt\vbox{ \hbox{\drwng{#1}{#2}{#3}}} }}$  }}  }
\def\rdet#1#2#3{\hbox{\phantom{}\kern-4pt\hbox{$\vcenter{\offinterlineskip \hbox
{
               \raise 4.5pt\hbox{\vrule height0.4pt width13pt depth0pt}
                \kern-1pt\vbox{ \hbox{\drwngt{#1}{#2}{#3}}} }}$  }}  }
 
\def\lde#1#2#3{\hbox{$\vcenter{\offinterlineskip  \hbox{
               \dlwng{#1}{#2}{#3}\kern-4.2pt\lower0.4pt\hbox{$\vcenter{\hrule 
                               width13pt}$}
               \kern-8pt\phantom{}   }}  $}}

\def\rwngb{\hbox{$\vbox{\offinterlineskip{
  \hbox{\phantom{}\kern6pt{$\bullet$}}\kern-2.5pt\hbox{$\Biggr/$}\kern-0.5pt
  \hbox{\phantom{}\kern-5pt$\bullet$}\kern-3.0pt\hbox{$\Biggr\backslash$}
  \kern-1.5pt\hbox{\phantom{}\kern6pt{$\bullet$}} }}$}}
 
\def\lwngb{\hbox{$\vbox{\offinterlineskip{ \hbox{$\bullet$}
  \kern-3.0pt\hbox{\phantom{}\kern6.0pt{$\Biggr\backslash$}}
  \kern-0.5pt\hbox{\phantom{}\kern11pt{$\bullet$}}\kern-3.5pt
  \hbox{\phantom{}\kern5.0pt {$\Biggr/$}}\kern-1.0pt\hbox{$\bullet$} }}$}}

\def\dbrwng#1#2#3{\hbox{$\vcenter{ \offinterlineskip{
  \hbox{\phantom{}\kern6pt{$\bullet^{\elbl{#3}}$}}
  \kern-2.5pt\hbox{$\Biggr/$}\kern-0.5pt
  \hbox{\phantom{}\kern-5pt$\bullet^{ \elbl{#1}}$}
  \kern-3.0pt\hbox{$\Biggr\backslash$}
  \kern-1.5pt\hbox{\phantom{}\kern6pt{$\bullet^{\elbl{#2}}$}}  } }$}}
 
\def\dblwng#1#2#3{\hbox{$\vcenter{\offinterlineskip{ \hbox{$\lelbl{#1}\bullet$}
  \kern-3.0pt\hbox{\phantom{}\kern6.0pt{$\Biggr\backslash$}}
  \kern-0.5pt\hbox{\phantom{}\kern11pt{$\lelbl{#2}\bullet$}}\kern-3.5pt
  \hbox{\phantom{}\kern5.0pt {$\Biggr/$}}\kern-1.0pt\hbox{$\lelbl{#3}\bullet$}}}
$} }

\def\rbde#1#2#3{\hbox{\phantom{}\kern-4pt\hbox{$\vcenter{\offinterlineskip  \hbo
x{
               \raise 4.5pt\hbox{\vrule height0.4pt width13pt depth0pt}
                \kern-1pt\vbox{ \hbox{\dbrwng{#1}{#2}{#3}}} }}$  }}  }
 
\def\lbde#1#2#3{\hbox{$\vcenter{\offinterlineskip  \hbox{
               \dblwng{#1}{#2}{#3}\kern-4.2pt\lower0.4pt\hbox{$\vcenter{\hrule w
idth13pt}$}
               \kern-8pt\phantom{}   }}  $}}

\def\ddgu#1.#2.{\dynk  \whtu0300{#1}\blku3000{#2}}
\def\ddgd#1.#2.{\dynk  \whtd0300{#1}\blkd3000{#2}}
 
\def\eddgiu#1.#2.#3.{\dynk \whtu0100{#1}\whtu1300{#2}\blku3000{#3}}
\def\eddgid#1.#2.#3.{\dynk \whtd0100{#1}\whtd1300{#2}\blkd3000{#3}}
\def\eddgiiu#1.#2.#3.{\dynk  \whtu0300{#1}\blku3100{#2}\blku1000{#3}}
\def\eddgiid#1.#2.#3.{\dynk  \whtd0300{#1}\blkd3100{#2}\blkd1000{#3}}

\def\ddfu#1.#2.#3.#4.{\dynk \whtu0100{#1}\whtu1200{#2}\blku2100{#3}\blku1000{#4}
}
\def\ddfd#1.#2.#3.#4.{\dynk \whtd0100{#1}\whtd1200{#2}\blkd2100{#3}\blkd1000{#4}
}
 
\def\eddfiu#1.#2.#3.#4.#5.{\dynk \whtu0100{#1}\whtu1100{#2}\whtu1200{#3}\blku210
0{#4}\blku1000{#5}}
\def\eddfid#1.#2.#3.#4.#5.{\dynk \whtd0100{#1}\whtd1100{#2}\whtd1200{#3}\blkd210
0{#4}\blkd1000{#5}}
\def\eddfiiu#1.#2.#3.#4.#5.{\dynk \whtu0100{#1}\whtu1200{#2}\blku2100{#3}\blku11
00{#4}\blku1000{#5}}
\def\eddfiid#1.#2.#3.#4.#5.{\dynk \whtd0100{#1}\whtd1200{#2}\blkd2100{#3}\blkd11
00{#4}\blkd1000{#5}}

\def\ddanu#1.#2.#3.#4.#5.{\dynk \whtu0100{#1}\whtu1100{#2}\whtu1100{#3}\cdots
                           \whtu1100{#4}\whtu1000{#5}}
\def\ddand#1.#2.#3.#4.#5.{\dynk \whtd0100{#1}\whtd1100{#2}\whtd1100{#3}\cdots
                           \whtd1100{#4}\whtd1000{#5}}
\def\ddandte#1.#2.#3.#4.#5.{\dynk \blkd0100{#1}\whtd1100{#2}\blkd1100{#3}\cdots
                           \blkd1100{#4}\whtd1000{#5}}
\def\ddandto#1.#2.#3.#4.#5.{\dynk \blkd0100{#1}\whtd1100{#2}\blkd1100{#3}\cdots
                           \whtd1100{#4}\blkd1000{#5}}
 
\def\eddanu#1.#2.#3.#4.#5.{\dynk \whtu0100{#1}\whtu1100{#2}%
                           \up1{\whtr0000{#3}}\cdots\whtu1100{#4}\whtu1000{#5}}
\def\eddand#1.#2.#3.#4.#5.{\dynk \whtd0100{#1}\whtd1100{#2}%
                           \up1{\whtr0000{#3}}\cdots\whtd1100{#4}\whtd1000{#5}}
\def\eddaid#1.#2.{\dynk\whtd0400{#1}\hskip30pt\whtd4000{#2}}

\def\eddanid#1.#2.#3.#4.#5.{\dynk \whtd0200{#1}\whtd2100{#2}%
                           \whtd1100{#3}\cdots\whtd1200{#4}\blkd2000{#5}}
\def\eddaniu#1.#2.#3.#4.#5.{\dynk \whtu0200{#1}\whtu2100{#2}%
                           \whtu1100{#3}\cdots\whtu1200{#4}\blku2000{#5}}

\def\eddaniid#1.#2.#3.#4.#5.#6.{\hbox{$\vcenter{\hbox
         {\dynk\hbox{$ \lbde{#1}{#2}{#3}\blkd1100{#4}\cdots%
          \blkd1200{#5}\whtd2000{#6} $}} }$}}
\def\eddaniiu#1.#2.#3.#4.#5.#6.{\hbox{$\vcenter{\hbox
         {\dynk\hbox{$ \lbde{#1}{#2}{#3}\blku1100{#4}\cdots%
          \blku1200{#5}\whtu2000{#6} $}} }$}}

\def\eddaiiid#1.#2.{\dynk\blkd0400{#1}\hskip30pt\whtd4000{#2}}

\def\ddbnu#1.#2.#3.#4.#5.{\dynk \whtu0100{#1}\whtu1100{#2}\whtu1100{#3}\cdots
                           \whtu1200{#4}\blku2000{#5}}
\def\ddbnd#1.#2.#3.#4.#5.{\dynk \whtd0100{#1}\whtd1100{#2}\whtd1100{#3}\cdots
                           \whtd1200{#4}\blkd2000{#5}}
 
\def\eddbnu#1.#2.#3.#4.#5.#6.{\dynk \lde{#1}{#2}{#3}\whtu1100{#4}\cdots
                           \whtu1200{#5}\blku2000{#6}}
\def\eddbnd#1.#2.#3.#4.#5.#6.{\dynk \lde{#1}{#2}{#3}\whtd1100{#4}\cdots
                           \whtd1200{#5}\blkd2000{#6}}

\def\ddcnu#1.#2.#3.#4.#5.{\dynk \blku0100{#1}\blku1100{#2}\blku1100{#3}\cdots
                           \blku1200{#4}\whtu2000{#5}}
\def\ddcnd#1.#2.#3.#4.#5.{\dynk \blkd0100{#1}\blkd1100{#2}\blkd1100{#3}\cdots
                           \blkd1200{#4}\whtd2000{#5}}
 
\def\eddcnu#1.#2.#3.#4.#5.#6.{\dynk \whtu0200{#1}\blku2100{#2}\blku1100{#3}
       \blku1100{#4}\cdots
                           \blku1200{#5}\whtu2000{#6}}
\def\eddcnd#1.#2.#3.#4.#5.{\dynk \whtd0200{#1}\blkd2100{#2}\blkd1100{#3}
       \cdots \blkd1200{#4}\whtd2000{#5}}

\def\dddnu#1.#2.#3.#4.#5.#6.{\hbox{$\vcenter{\hbox
         {\dynk\hbox{$ \whtu0100{#1}\whtu1100{#2}\cdots%
          \whtu1100{#3}\rde{#4}{#5}{#6} $}}  }$}}
\def\dddnd#1.#2.#3.#4.#5.#6.{\hbox{$\vcenter{\hbox
         {\dynk\hbox{$ \whtd0100{#1}\whtd1100{#2}\cdots%
          \whtd1100{#3}\rde{#4}{#5}{#6} $}} }$}}
\def\dddndte#1.#2.#3.#4.#5.#6.{\hbox{$\vcenter{\hbox
         {\dynk\hbox{$ \blkd0100{#1}\whtd1100{#2}\cdots%
          \blkd1100{#3}\rdet{#4}{#5}{#6} $}} }$}}
\def\dddndto#1.#2.#3.#4.#5.#6.{\hbox{$\vcenter{\hbox
         {\dynk\hbox{$ \whtd0100{#1}\blkd1100{#2}\cdots%
          \blkd1100{#3}\rdet{#4}{#5}{#6} $}} }$}}
\def\dddiv#1.#2.#3.#4.{\hbox{$\vcenter{\hbox
         {\dynk\hbox{$ \whtu0100{#1}\rde{#2}{#3}{#4}
              $}}  }$}}

\def\edddnu#1.#2.#3.#4.#5.#6.#7.#8.{\hbox{$\vcenter{\hbox
         {\dynk\hbox{$ \lde{#1}{#2}{#3}\whtu1100{#4}\cdots%
          \whtu1100{#5}\rde{#6}{#7}{#8} $}}  }$}}
\def\edddnd#1.#2.#3.#4.#5.#6.#7.#8.{\hbox{$\vcenter{\hbox
         {\dynk\hbox{$ \lde{#1}{#2}{#3}\whtd1100{#4}\cdots%
          \whtd1100{#5}\rde{#6}{#7}{#8} $}} }$}}

\def\edddniid#1.#2.#3.#4.#5.{\hbox{$\vcenter{\hbox
         {\dynk\hbox{$ \blkd0200{#1}\whtd2100{#2}\whtd1100{#3}\cdots%
          \whtd1200{#4}\blkd2000{#5} $}} }$}}
\def\edddniiu#1.#2.#3.#4.#5.{\hbox{$\vcenter{\hbox
         {\dynk\hbox{$ \blku0200{#1}\whtu2100{#2}\whtu1100{#3}\cdots%
          \whtu1200{#4}\blku2000{#5} $}} }$}}

\def\ddei#1.#2.#3.#4.#5.#6.{\hbox{$\vcenter{\hbox
       {\dynk \whtd0100{#1}\whtd1100{#3}%
       \up1{\whtr0001{#2}}\whtd1110{#4}\whtd1100{#5}\whtd1000{#6}} }$}}
\def\ddeit#1.#2.#3.#4.#5.#6.{\hbox{$\vcenter{\hbox
       {\dynk \whtd0100{#1}\blkd1100{#3}%
       \up1{\blkr0001{#2}}\whtd1110{#4}\blkd1100{#5}\whtd1000{#6}} }$}}
 
\def\eddei#1.#2.#3.#4.#5.#6.#7.{\hbox{$\vcenter{\hbox
       {\dynk \whtd0100{#1}\whtd1100{#3}%
       \up1{\whtr0011{#2}}\up2{\whtr0001{#7}}\whtd1110{#4}\whtd1100{#5}%
       \whtd1000{#6}} }$}}

\def\ddeii#1.#2.#3.#4.#5.#6.#7.{\hbox{$\vcenter{\hbox
       {\dynk \whtd0100{#1}\whtd1100{#3}%
       \up1{\whtr0001{#2}}\whtd1110{#4}\whtd1100{#5}\whtd1100{#6}%
       \whtd1000{#7}} }$}}
\def\ddeiit#1.#2.#3.#4.#5.#6.#7.{\hbox{$\vcenter{\hbox
       {\dynk \whtd0100{#1}\blkd1100{#3}%
       \up1{\blkr0001{#2}}\whtd1110{#4}\blkd1100{#5}\whtd1100{#6}%
       \blkd1000{#7}} }$}}
 
\def\eddeii#1.#2.#3.#4.#5.#6.#7.#8.{\hbox{$\vcenter{\hbox
       {\dynk \whtd0100{#8}\whtd1100{#1}\whtd1100{#3}%
       \up1{\whtr0001{#2}}\whtd1110{#4}\whtd1100{#5}\whtd1100{#6}%
       \whtd1000{#7}} }$}}

\def\ddeiii#1.#2.#3.#4.#5.#6.#7.#8.{\hbox{$\vcenter{\hbox
       {\dynk \whtd0100{#1}\whtd1100{#3}%
       \up1{\whtr0001{#2}}\whtd1110{#4}\whtd1100{#5}\whtd1100{#6}%
       \whtd1100{#7}\whtd1000{#8}} }$}}
\def\ddeiiit#1.#2.#3.#4.#5.#6.#7.#8.{\hbox{$\vcenter{\hbox
       {\dynk \whtd0100{#1}\blkd1100{#3}%
       \up1{\blkr0001{#2}}\whtd1110{#4}\blkd1100{#5}\whtd1100{#6}%
       \blkd1100{#7}\whtd1000{#8}} }$}}
 
\def\eddeiii#1.#2.#3.#4.#5.#6.#7.#8.#9.{\hbox{$\vcenter{\hbox
       {\dynk \whtd0100{#1}\whtd1100{#3}%
       \up1{\whtr0001{#2}}\whtd1110{#4}\whtd1100{#5}\whtd1100{#6}%
       \whtd1100{#7}\whtd1100{#8}\whtd1000{#9}} }$}}

\def\bb{{\cal B}}
\def\ff{{\cal F}}\def\gg{{\cal G}}
\def\ii{{\cal I}}
\def\mm{{\cal M}}\def\nn{{\cal N}}\def\oo{{\cal O}}
\def\ss{{\cal S}}                                     
\def\uu{{\cal U}}

\def\c{{\mathbb C}}
\def\h{{\mathbb H}}
\def\o{{\mathbb O}}\def\bp{{\mathbb P}}
\def\r{{\mathbb R}}\def\z{{\mathbb Z}}

\def\gothg{{\mathfrak g}}\def\gothh{{\mathfrak h}}
\def\so{{\mathfrak so}}

\def\one{{\bf 1}}\def\two{{\bf 2}}\def\three{{\bf 3}}\def\four{{\bf 4}}
\def\six{{\bf 6}}\def\seven{{\bf 7}}
\def\bolde{{\bf e}}\def\boldf{{\bf f}}
\def\boldl{{\bf l}}\def\boldn{{\bf n}}\def\boldx{{\bf x}}

\def\ftil{{\widetilde f}}
\def\ks{{\rm ks}}
\def\sch{{\rm Sch}}
\def\homo{{\rm Hom}}
\def\endo{{\rm End}\ }
\def\aut{{\rm Aut}\ }
\def\pic{{\rm Pic}}
\def\nm{{\rm Nm}}
\def\spin{{\rm Spin}}
\def\grass{{\rm Grass}}
\def\gr{{\rm gr}}
\def\fix{{\rm Fix}}
\def\gal{{\cal G}}
\def\ad{{\rm ad}\ }
\def\sdp{\mathbin{\hbox{\twelvesymb o}}}
\def\idfy{\mathbin{\hbox{$\,\widetilde\rightarrow\,$}}}
\def\dual{^{\vee}}

\def\im{{\rm Im\ }}
\def\mat{{\rm Mat}}

\def\ker{{\rm ker \ }}
\def\rk{{\rm rank \ }}
\def\tr{{\rm trace}\ }

\def\peta{P_{\eta}}
\def\su{{\cal SU}_C}
\def\map#1{\ \smash{\mathop{\longrightarrow}\limits^{#1}}\ }

\def\ses#1#2#3{0\rightarrow {#1}\rightarrow {#2}\rightarrow {#3}\rightarrow 0}
\def\prf{{\it Proof.}\ }
\def\deq{\mathrel{\lower3pt\hbox{$\buildrel{\rm def}\over=$}}}

\def\qqed{{\unskip\nobreak\hfill\hbox{ $\Box$}\medskip\par}}

\def\<{\langle}
\def\>{\rangle}

\begin{document}

\title{Trigonal curves and Galois $\spin(8)$-bundles}
\author{W.M. Oxbury and S. Ramanan}
\date{}

\maketitle

\bigskip

The group $\spin(8)$ occupies a special position among the complex simple Lie 
groups, in having outer automorphism group $S_3$ and an appealing 
Dynkin diagram: 
$$
\dynk\dddiv{V}..{S^-}.{S^+}.
$$
The outer nodes are dual to three fundamental 
8-dimensional representations, the standard representation $V$ on 
which the group acts via the double cover $\spin(8) \rightarrow SO(8)$, 
and the half-spinor representations $S^{\pm}$. These spaces are 
permuted by the action of $S_3$, and between any two 
there is a Clifford multiplication to the third. It 
is possible to identify all three spaces with the complex Cayley 
algebra $\o$, and Clifford multiplication with multiplication of 
octonions. In this sense $\spin(8)$ is a member of the `exceptional' 
club, and is closely related to the other exceptional groups 
$G_2 = \aut \o,F_4,E_6$ etc.

\medskip

In this paper we study the moduli of principal holomorphic $\spin(8)$-bundles
over an algebraic curve. However, in order to exploit triality, and to obtain 
a moduli space with particularly nice properties,  
we impose some additional constraints on our bundles. 
We suppose that $X$ is an algebraic curve with an $S_3$-action. The 
group $S_3$ 
then acts in two ways on $\spin(8)$-bundles over a curve: by 
pull-back 
under the action on $X$, and by the triality action on the structure 
group. We call a $\spin(8)$-bundle {\it Galois} if these two actions coincide,
that is, if it is a fixed point of the group action
$$
F\mapsto u^*F^u,
\qquad
u\in S_3,
$$
on the moduli variety $\mm_X(\spin(8))$ of (semistable) $\spin(8)$-bundles. 
We denote the fixed-point set by $\ff_X \subset \mm_X(\spin(8))$, and 
distinguish the subvariety of $\ff_X$ consisting of bundles which admit a 
lift of the $S_3$-action (see Lemma \ref{galoisbundle} for the precise notion). 
This subvariety has a partial desingularisation $\nn_X$
parametrising pairs $(F,\Lambda)$ where $F\in \ff_X$ and $\Lambda$ is an
$S_3$-lift, or more precisely a splitting $\Lambda: S_3 \rightarrow \gg_F$
of the Mumford group associated to $F$ (see Definition \ref{definition}).
$\nn_X$ is the main object of interest of this paper.

\medskip

In \S\ref{moduli} we compute the tangent spaces at stable points and the 
semistable boundary of $\nn_X$ (see \ref{tgt-space} and \ref{ssboundary}).
The best situation, to which we mainly restrict ourselves, is where the quotient
$X/S_3$ is isomorphic to $\bp^1$, and in this case we show that $\nn_X$ is 
smooth at stable Galois bundles and that the locus of 
nonstable bundles can be identified with the moduli
space $\ss\uu_{X/\sigma}(2)$ of rank 2, trivial determinant vector bundles
on the quotient of $X$ by the involution $\sigma \in S_3$ which exchanges 
the spinor representations (see \ref{notation} for notation). 
Both properties are in marked contrast to 
$\mm_X(\spin(8))$ itself, whose semistable boundary is more 
complicated to describe
and whose stable points are not necessarily smooth but may have 
finite quotient singularities.

\medskip

The condition $X/S_3 \cong \bp^1$ just means that $X = \gal (C)$ 
is the Galois closure over $\bp^1$
of the trigonal curve $C = X/\sigma$, and conversely 
one can construct $\nn_X$ starting from any 
(non-cyclic) trigonal curve $C\map{3:1} \bp^1$. Indeed, this was 
precisely our motivation 
for the construction, and we shift attention to the trigonal point of view in 
\S\ref{galois3}. The curve
$X = \gal (C)$ is a branched double cover of $C$; we write
$\nn_C = \nn_X$; and the main result of the paper is the following:

\statement Theorem \ref{dimtheorem}. Given a (non-cyclic) trigonal curve 
$C\map{3:1} 
\bp^1$ there exists a projective moduli space $\nn_C$ 
(parametrising Galois $\spin(8)$-bundles on $\gal(C)$)
which admits an 
inclusion $\su(2) \hookrightarrow \nn_C$ such that $\nn_C$ is smooth and of 
dimension $7g-14$ away from the image of $\su(2)$. 

(Cyclic trigonal curves need to be treated separately as the Galois group
is $\z/3$ rather than $S_3$, and although many of the computations work well 
here---and we expect that the above theorem remains true---we 
confine ourselves to brief remarks on this case at various points---see 
the end of \S\ref{hypcase}.)

\medskip

The moduli space $\su(2)$ rank 2 semistable vector bundles with trivial 
determinant is an object of considerable interest. 
As well as containing the Jacobian Kummer as its singular locus, 
a fundamental geometric feature 
is the so-called {\it Schottky configuration} of Prym Kummers.
(See for example \cite{vGP}, \cite{NR3}.)
That is, the group $J_C[2]$ of 2-torsion points in the Jacobian acts on
$\su(2)$ by tensor product; the fixed-point set of an element $\eta \in J_C[2]$
is a pair of (isomorphic)
Kummer varieties $(P_{\eta}/\pm) \cup (P^-_{\eta}/\pm) $ of dimension
$g-1$. $P_{\eta}$ is the Prym variety of the unramified double cover 
$C_{\eta} \rightarrow C$ corresponding to $\eta$, and the map 
$P_{\eta} \cup P^-_{\eta}\rightarrow \su(2)$
is (up to a choice of $\eta^{1\over 2}$) the 
direct image of line bundles from $C_{\eta}$. The incidence
relations among these Kummer varieties as $\eta$ varies, interpreted via 
the embedding $\su(2) \hookrightarrow |2\Theta|$ (where $\Theta$ is the Riemann
theta divisor in the Jacobian) correspond precisely to the Schottky-Jung-Donagi
identities among the thetanulls of $J_C$ and the Pryms.

\medskip

It turns out that when $C$ is trigonal 
the Schottky configuration too has a `fattening' in
the moduli space $\nn_C$. By a beautiful and well-known construction of 
Recillas \cite{Rec}, each trigonal Prym $P_{\eta}$ is isomorphic as a 
principally polarised abelian variety to the Jacobian of a tetragonal 
curve $R_{\eta} \map{4:1} \bp^1$. Thus a trigonal Schottky configuration
consists of Jacobian Kummers $J_{R_{\eta}}/\pm$. Each of these 
is also the singular locus of a moduli space of bundles $\ss\uu_{R_{\eta}}(2)$.
We shall show: 

\statement Theorem \ref{nonab-schottky}. Given a trigonal curve $C$ and 
a nonzero 2-torsion point $\eta \in J_C[2]$ there exists 
(up to a choice of $\eta^{1\over 2}$) a natural map $\ss\uu_{R_{\eta}}(2)
\hookrightarrow\nn_C$
for which the following diagram commutes:
$$
\begin{array}{ccc}
P_{\eta}& \idfy J_{R_{\eta}} \rightarrow &\ss\uu_{R_{\eta}}(2)\\
&&\\
\downarrow && \downarrow\\
&&\\
\ss\uu_C(2)&\hookrightarrow& \nn_C.\\
\end{array}
$$

One can view the right-hand side of this diagram, as $\eta$ varies, 
as a nonabelian
Schottky configuration singular along the classical Schottky configuration 
on the left-hand side.

\medskip

Finally, we wish to add a word about the motivation for these constructions. 
Originally, this was our interest 
the projective embedding (see \cite{BV}, \cite{vGI}; 
we assume $C$ is 
nonhyperelliptic unless $g=2$)
$$
\phi: \su(2) \rightarrow |2\Theta| = \bp^{2^g -1}.
$$
It is well-known that for $g=2$ the map $\phi$ is an isomorphism 
$\su(2) \idfy \bp^3$ (see \cite{NR1});
and that for $g=3$ (see \cite{NR2}) the image $\phi(\su(2)) \subset \bp^7$ 
is the unique 
Heisenberg-invariant quartic---the  {\it Coble quartic}---singular 
along the Kummer variety $J_C/\pm \subset 
|2\Theta|$. When $g=4$ it was shown in \cite{OP} that in $|2\Theta| = \bp^{15}$
there exists a unique Heisenberg-invariant quartic singular along 
the image $\phi(\su(2))$. It was shown, moreover, that the fixed-point set of
$\eta \in J_C[2]$ in this quartic is a pair of Coble quartics 
$\ss\uu_{R_{\eta}}(2)$ of the corresponding Recillas curve. However, the 
questions remained open: what is this 14-dimensional quartic as a moduli space?,
and is its singular locus equal to $\su(2)$?

\medskip

The crucial point
here is that a general curve of genus 4 is trigonal (in two ways. 
However, by Torelli's theorem $R_{\eta}$ is independent of the choice of 
trigonal structure). In view of the results of the present paper
it is natural to expect that the genus 4 quartic in $\bp^{15}$ is $\nn_C$.
We hope to pursue this question in a sequel.

\medskip
\noindent
{\it Acknowledgments:} 
The first author wishes to thank TIFR, Mumbai for its hospitality 
during a visit in July 1998 when much of this work was carried out.
We are also grateful to Steve Wilson for pointing out Lemma 
\ref{wilson}.
In writing this paper the authors were partially supported by EPSRC
grants GR/M03924 and GR/M36663.

\section{$\spin(8)$ and triality}
\label{spin}

We begin by recalling the triality story for $\spin(8)$ from various points of
view. Roughly, \S\ref{group} will be used in the discussion of 
Galois bundles of section \ref{trigonal}; 
\S\ref{2by2} and \S\ref{6quadrics} are needed for the 
discussion of stability in section \ref{moduli}; and \S\ref{liealg} will be 
used in the 
computation of the dimension of the moduli space in section \ref{galois3}.

\subsection{The group $\spin(8)$}
\label{group}

We shall always denote by 
$V \cong \c^{2n}$ the standard orthogonal representation of 
$\spin({2n})$, and by $S^{\pm} \cong \c^{2^{n-1}}$ its half-spinor 
representations. We shall denote by $q$ the quadratic form on $V$; 
if we fix a decomposition $V = N \oplus N\dual$, where $N$ and 
$N\dual$ are maximal isotropic subspaces, dual via $q$, then by definition
\begin{equation}
\label{spinors}
\textstyle 
S^{+} = \bigwedge^{\rm even}N,
\qquad
S^{-} = \bigwedge^{\rm odd}N.
\end{equation}

We need to recall some facts about 
$S^{\pm}$
(see \cite{C}).
Start with the bilinear pairing
$$
\begin{array}{rcl}
	r:\bigwedge N \otimes \bigwedge N & \rightarrow & \textstyle 
\bigwedge^{n}N 
	\cong \c  \\
	s\otimes t & \mapsto & (\beta(s) \wedge t)_{n}  \\
\end{array}
$$
where $\beta$ is the principal anti-involution of $\bigwedge N$, i.e. 
it is the identity on $N$ and reverses multiplication; and where 
$()_{n}$ denotes the component of top degree. When $n$ is odd this 
form vanishes on each of $S^{\pm}$ and induces a nondegenerate 
pairing $r:S^{+} \otimes S^{-} \rightarrow \c$. When $n$ is even it 
restricts to a nondegenerate pairing on each of $S^{\pm}$, which is 
symmetric precisely when $n\equiv 0(4)$. 

We shall be 
concerned with the case $n=4$; 
here the pairing $r$ determines quadratic forms $q^{\pm}$ on
the spinor spaces $S^{\pm}$. Moreover, these determine an embedding
(in fact the projection into the product of any two factors is still 
an embedding)
\begin{equation}
\label{spin8subgp}
\spin(8) \hookrightarrow SO(V) \times SO(S^{+}) \times SO(S^{-}).
\end{equation}

There is a well-known triality relationship among
the three orthogonal factors here. 
$\spin(8)$ has outer automorphism group $S_{3}$, and
we shall recall explicitly how $S_3$ acts on $\spin(8)$, or more precisely 
how to split the sequence
$$
1\rightarrow {\rm Inn}\ \spin(8) \rightarrow {\rm Aut}\ \spin(8) 
\rightarrow S_3
\rightarrow 1.
$$

\begin{notat}\rm
\label{notation}
We shall use, here and in later sections,
generators $\sigma = (23)$ and $\tau = (123) \in S_{3}$; 
thus $\sigma$ should interchange $S^{\pm}$, as an involution of 
$\spin(8)$; while $\sigma \tau: V \leftrightarrow S^{+}$ and 
$\sigma \tau^{2}: V \leftrightarrow S^{-}$. 
\end{notat}

Choose unit vectors
$v_{0} \in V$, $s_{0}\in S^{+}$ and $t_{0}= v_{0} \cdot s_{0}\in  
S^{-}$. 
Given these choices we let $S_3$ act on the vector space $V\oplus 
S^{+}\oplus S^{-}$ as follows. 
$\sigma$ acts by reflection of $V$ in $v_{0}^{\perp}$ and by Clifford 
multiplication $v_{0}: S^{+}\leftrightarrow S^{-}$. Similarly the 
group elements $\sigma \tau$ and $\sigma \tau^{2}$ are represented by 
$s_{0}$ and $t_{0}$ respectively. 

Given these choices we shall denote the resulting representation by 
$$
\rho : S_3 \rightarrow O(V\oplus S^+ \oplus S^-).
$$
This intertwines an $S_3$-action on 
$\spin(8)$ by $g\mapsto g^u := \rho (u)^{-1}g
\rho(u)$ for $u\in S_3$ and $g \in 
\spin(8)\subset O(V\oplus S^+ \oplus S^-)$. Note that the subgroup 
of $O(V\oplus S^+ \oplus S^-)$ generated by the images of $S_3$ and 
$\spin(8)$ is isomorphic to the semi-direct product $\spin(8) \sdp 
S_3$. 
We note for later use that the group law on $\spin(8) \sdp S_3$ is
\begin{equation}
\label{sdp}
(g,u)(h,v)= (gh^u,uv)
\qquad
u,v \in S_3
\quad
g,h \in \spin(8).
\end{equation}

In view of triality, there is a 
multiplication map $S^{+}\otimes S^{-} \rightarrow V$
permuted by $S_3$ with the Clifford multiplications 
$V\otimes S^{\pm} \rightarrow S^{\mp}$.
This can be
described in terms of a trilinear form:
\begin{equation}
	\label{trilinear}
	\begin{array}{rcl}
		c: V\otimes S^{+}\otimes S^{-} & \rightarrow & \c  \\
		v\otimes s\otimes t & \mapsto & q^{-}(v\cdot s,t) = q^{+}(v\cdot 
t,s) 
		\deq q(s\cdot t,v). 
	\end{array}
\end{equation}
The first equation is easily checked; the second defines the 
Clifford multiplication by:
\begin{displaymath}
	\begin{array}{ccc}
		S^{+}\otimes S^{-} & \map{c} & V\dual  \\
		&&\\
		 &.\searrow & \downarrow q  \\
		 &&\\
		 &  & V
	\end{array}
\end{displaymath}

\begin{rem}\rm
One can use the three Clifford multiplications to define a 
commutative algebra 
structure on $V\oplus S^+ \oplus S^-$ by taking the product of two 
vectors 
in the same summand to be zero. The resulting algebra, called the 
{\it Chevalley algebra}, is an example of a vertex operator algebra. 
(See 
\cite{FFH}.)
\end{rem}

\begin{lemm}
\label{outerautos}
The orthogonal action $\rho: S_{3} \rightarrow O(V\oplus 
S^{+}\oplus S^{-})$ 
preserves the cubic form induced by the 
trilinear form $c$, and hence also the Chevalley algebra structure, 
and determines an $S_3$-action on 
$\spin(8)$ as above. Moreover, every splitting $S_3 \rightarrow {\rm
Aut}\ 
\spin(8)$ arises in this way.
\end{lemm}

\begin{rem}\rm
\label{oct1}
Given the above choices $V$ acquires an algebra structure
$$
\begin{array}{rcl}
	V\times V & \rightarrow & S^{+}\times S^{-} \map{\cdot}  V  \\
	(u,v) & \mapsto  & (t_{0}\cdot u,s_{0}\cdot v) 
\end{array}
$$
It is well-known 
(and will follow from Remark \ref{oct2} below) that
this is precisely the (complex) Cayley algebra $V\cong \o$ with centre 
$\<v_{0}\>$. 
One deduces from this fact a characterisation of the image of
(\ref{spin8subgp}): 
\begin{equation}
\label{spin8subgp2}
\spin(8) = \{ (a,b,c)\in SO(\o) \times SO(\o) \times SO(\o) \ |\ 
a(u)b(v) = c(uv)
\quad
\forall \ u,v\in \o \}.
\end{equation}
In this language the triality action is given by $\sigma: 
(a,b,c)\mapsto (a',c',b')$ and $\tau: (a,b,c)\mapsto (b',c,a')$ where 
$a'(u) = \overline{a(\overline u)}$.  
\end{rem}

\begin{lemm}
\label{zero-norm}
For an isotropic subspace $U\subset V$ and an involution $g\in S_3$
not preserving $V$, consider the Clifford 
multiplication maps:
$$
\mu: U\otimes U \map{{\rm id}\otimes g} V \otimes S^{\pm}
\map{\cdot} S^{\mp}
$$
If $\dim U =1$ then $\mu =0$, and if $\dim U =2$ then $\rk \mu \leq 1$.
\end{lemm}

\prf
Suppose first that $\dim U =1$. 
Identifying $V,S^+,S^-$ with $\o$ as in Remark \ref{oct1} the map $\mu$ 
becomes $x\otimes x\mapsto x\overline x = \| x\| =0$. (See also Remark
\ref{oct2}). If $\dim U =2$ then the symmetric tensors $S^2 U \subset 
U\otimes U$ are spanned by squares, and so are contained in kernel of
$\mu$. This shows that $\rk \mu \leq 1$.
\qqed

Finally, note that the centre of $\spin(8)$ is $Z(\spin(8)) = \z/2 
\times \z /2$, consisting of the matrices of 
$8\times 8$ blocks:
$$
\pmatrix{1&&\cr &1& \cr &&1\cr},
\pmatrix{1&&\cr &-1& \cr &&-1\cr},
\pmatrix{-1&&\cr &1& \cr &&-1\cr},
\pmatrix{-1&&\cr &-1& \cr &&1\cr}.
$$

\subsection{Triality in terms of $2\times 2$ matrices}
\label{2by2}

We can always split $V= \c^8$, as an orthogonal space,
into a direct sum of orthogonal 
$\c^4$s. In turn $\c^4$ can be identified with $\homo(A,B)$ 
where $A,B \cong \c^2$, with orthogonal structure given by the 
determinant (geometrically, every smooth complex quadric surface is 
the Segre $\bp^1 \times \bp^1$). The spaces $A\dual$ and $B$ are then 
the spinor representations, where $\spin(4) = SL(2) \times SL(2)$. 
We can view this as saying that $A,B$ carry fixed (complex) 
orientations
$\lambda_A: \bigwedge^2 A \idfy \c$ and 
$\lambda_B: \bigwedge^2 B \idfy \c$.

Let us generalise this situation for a moment. Suppose that $A,B 
\cong \c^n$ and are each equipped with a fixed orientation. Then any 
$u\in \homo(A,B)$ has an {\it adjugate} homomorphism $\overline u \in 
\homo(B,A)$. This is just the transpose of the composition
$$
\textstyle
A\dual \map{\lambda_A} \bigwedge^{n-1} A \map{\bigwedge^{n-1} u}
\bigwedge^{n-1} B \map{\lambda_B} B\dual,
$$
and with respect to chosen bases is given by the transpose matrix of 
signed cofactors of the matrix representing $u$. Adjugacy is 
therefore a natural birational 
involution $\homo(A,B) \leftrightarrow \homo(B,A)$
for oriented vector spaces; when $n=2$ it is a linear involution.
(The reader may care to consider the next case $n=3$: projectively 
the map blows up the Segre $\bp^2 \times \bp^2$ and contracts secant 
lines down to the dual Segre fourfold.)

In the case $n=2$ the determinant of a map between oriented spaces
can be interpreted as 
a quadratic form $\det: \homo(A,B) \rightarrow \c$ whose polarisation is the 
symmetric bilinear form
$$
\<u,v\> = {1\over 2} \tr u\overline v,
\qquad
u,v \in \homo(A,B).
$$

We now return to our 8-dimensional orthogonal space, which we shall 
decompose as 
$$
V= \homo(A,B) \oplus \homo(C,D)
$$
where $A,B,C,D$ are oriented $\c^2$s.
It follows from Lemma \ref{spin-sum}, and the identifications $A\idfy 
A\dual$ etc given by the orientations, that the spinor spaces are
$$
 \begin{array}{rcl}
 	S^+ & = &  \homo(B,C) \oplus \homo(A,D), \\
 	S^- & = & \homo(C,A)\oplus \homo(B,D).
 \end{array}
$$
It is straightforward to calculate the Clifford multiplication maps:
we shall take on $V$ the quadratic form $(a,b) \mapsto \det a - \det b$.
One then finds:
$$
\begin{array}{rcl}
	V\otimes S^+ & \rightarrow & S^-  \\
	(a,b)\otimes (x,y) & \mapsto  & 
	(\overline{a}\ \overline{x} +\overline{y} 
	b,  bx + y\overline{a})
\end{array}
$$
with similar expressions for the other two maps. (The rule is: for 
each of the two Hom summands, simply add the only possible 
composites, allowing adjugates, so that the expression is well-defined.)

\begin{rem}\rm
\label{rankdrops}
The orientations on the spaces $A,B,C,D$ induce identifications
$\bigwedge^8 V \cong \bigwedge^8 S^+ \cong \bigwedge^8 S^-$, so the 
determinant of each Clifford multiplication map is a well-defined scalar.
For each pair of $2\times 2$ matrices $(a,b) \in V = \homo(A,B)\oplus 
\homo(C,D)$ one 
finds that the linear map $m_{a,b} : S^+ \rightarrow 
S^-$ satisfies: 
$$
\det m_{a,b} = -(\det a - \det b)^4.
$$
This corresponds to the fact that the rank drops by 4 if
$\det a = \det b$, i.e. if $(a,b)\in V$ is isotropic. One can check this 
directly (e.g. using Maple) by Gauss-Jordan elimination.
\end{rem}

\begin{rem}\rm
\label{oct2}
If we fix unimodular isomorphisms $A\idfy B \idfy C \idfy D$ 
we obtain a natural choice 
of unit vectors $v_{0} = (1,0) \in V$, $s_{0} = (1,0) \in S^+$, 
$t_{0} = (1,0) \in S^-$. (Note that $t_0 = v_0 s_0$ under the Clifford 
multiplication above.) With this choice
the algebra structure of Remark \ref{oct1} becomes:
$$
\begin{array}{rcl}
	V\times V & \rightarrow & V  \\
	(a,b)(c,d) & \mapsto & (ac+\overline{d} b,da+b\overline{c}).
\end{array}
$$
On the other hand, the chosen isomorphisms identify $V = \endo A 
\oplus \endo A$ where $\endo A \cong \mat_{2}(\c) \cong 
\h \otimes_{\r}\c$, the complexification of the quaternions. Moreover, 
quaternionic conjugation is precisely the map which sends a matrix to 
its adjugate. 
It follows that the above algebra structure on $V$ is exactly the 
(split) Cayley-Dickson process applied to $\h \otimes_{\r}\c$ (see 
\cite{harvey} pp.105--106), i.e. identifies $V$ with $\o$.
\end{rem}

\subsection{6-dimensional quadrics}
\label{6quadrics}

Let $Q\subset 
V$ and $Q^{\pm} \subset S^{\pm}$ be the quadrics
defined by the respective quadratic 
forms on the basic representations of $\spin(8)$.
Each quadric has two families of 4-dimensional isotropic 
spaces, which we shall sometimes refer to as {\it $\alpha$-planes} 
and 
{\it $\beta$-planes}. 
Such subspaces $A,B \subset Q$ belong to 
opposite families if and only if
$$
\dim A\cap B \equiv 1
\quad
{\rm mod}\ 2.
$$
Moreover, these families are parametrised precisely by the quadrics $Q^{\pm}$
as follows.

Given vectors $s\in S^+$ and $t\in S^-$, consider the Clifford multiplication
maps $m_s: V \leftrightarrow S^-$ and $m_t: V \leftrightarrow S^+$. We have 
observed in Remark \ref{rankdrops} that these maps have nonzero kernel
exactly when $s,t$ respectively are isotropic, and that then the rank is 4. 
We define:
$$
\begin{array}{rcl}
A_s &=& \ker \{V\map{m_s} S^-\} = \im \{S^- \map{m_s} V\} \subset Q \subset V,\\
B_t &=& \ker \{V\map{m_t} S^+\} = \im \{S^+ \map{m_t} V\} \subset Q \subset V.\\
\end{array}
$$
One readily checks (see \cite{C}, or use the set-up of of \S\ref{2by2})
these equalities, that they define isotropic subspaces, and that:
$$
\dim A_s \cap B_t = 
\cases{1& spanned by $s\cdot t \in V$ if $ s\cdot t \not=0$,\cr
3& if $ s\cdot t =0$.\cr}
$$
In particular $A_s, B_t$ are in opposite families. 

By triality 
(i.e. using (\ref{trilinear}) in \S\ref{spin}) we can make 
the same constructions in each of $S^{\pm}$. We summarise our notation
in the following diagram:
\begin{equation}
\label{quad-pic}
 \leavevmode
\epsfxsize=3in\epsfbox{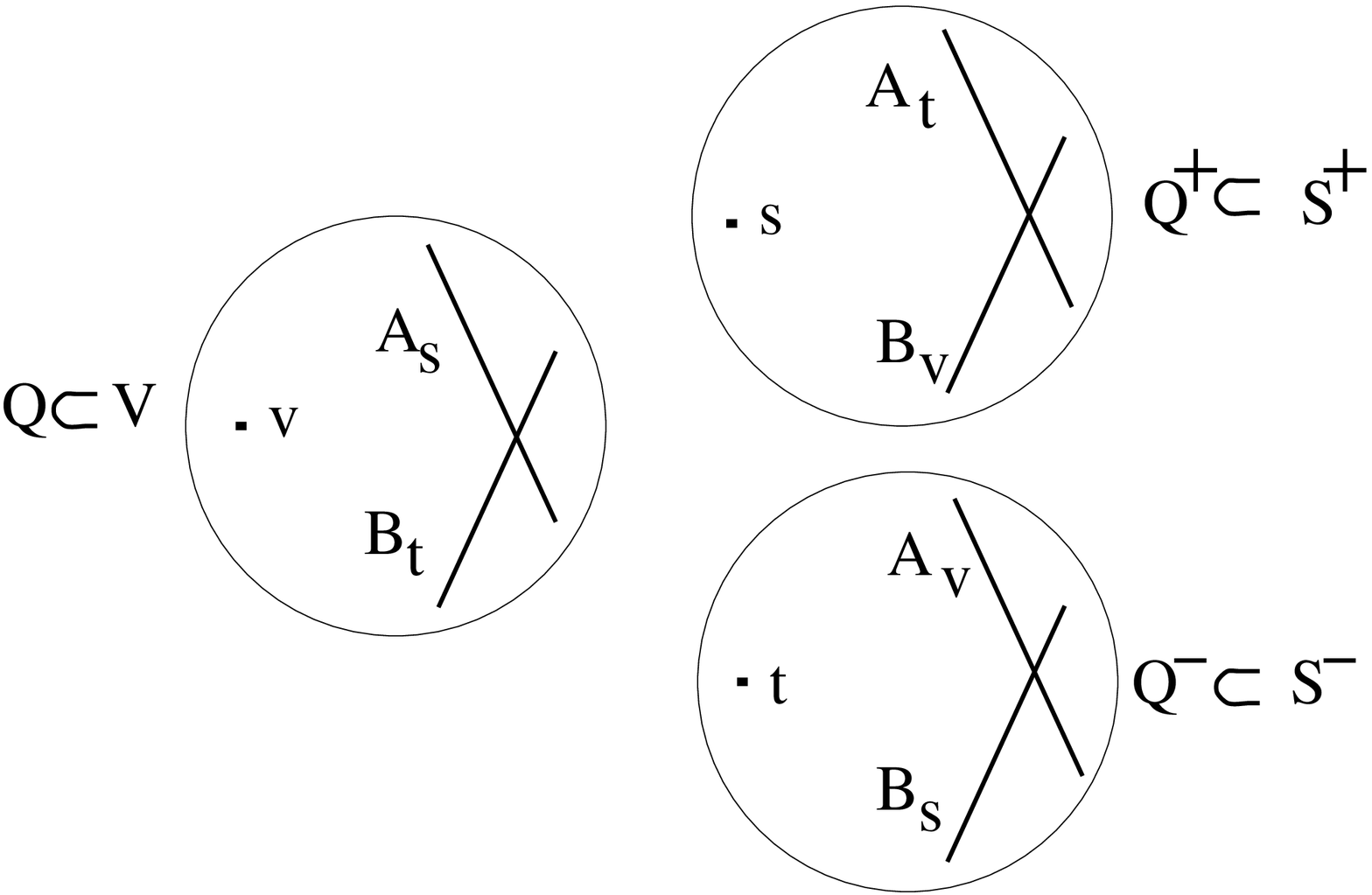}
\end{equation}
Note that:
\begin{equation}
\label{prod=0}
s\cdot t =0
\quad \Longleftrightarrow\quad 
s\in A_t
\quad \Longleftrightarrow\quad 
t\in B_s
\quad \Longleftrightarrow\quad 
\dim A_s \cap B_t = 3.
\end{equation}
Moreover, there are canonical dualities $A_s = B_s\dual$ for $s\in Q^\pm$ 
or $Q$; this follows, for example, via the quadratic form on $V$, from 
the exact sequences given by Clifford multiplication:
\begin{equation}
\label{ABsequences}
\begin{array}{c}
\ses{A_s}{V}{B_s},\\
\ses{B_t}{V}{A_t}.\\
\end{array}
\end{equation}
Geometrically, of course, $A_t$ is the set of $s\in Q^+$ for which 
$\dim A_s \cap B_t = 3$, and so parametrises hyperplanes in $B_t$.

Similarly, for distinct $s,s' \in Q^+$ we have
\begin{equation}
\label{polarity}
\dim A_s \cap A_{s'} =
\cases{2 & if $\<s,s'\> =0$,\cr
       0 & otherwise.\cr}
\end{equation} 
Note that this fact identifies the Grassmannian $\grass(2,A_s)$ with the 
4-dimensional (tangent) quadric in $\bp (s^\perp /s)$.

We can generalise the above correspondence by defining, 
for any isotropic subspace $U\subset V$, and using the 
trilinear form (\ref{trilinear}):
$$
\begin{array}{rcl}
	A_U = \ker\{S^- \map{c} (U\otimes S^+)\dual\} & =&
\im \{U\otimes S^+ \map{c} (S^-)\dual\}^{\perp},     \\
B_U = \ker\{S^+ \map{c} (U\otimes S^-)\dual\} & =&
\im \{U\otimes S^- \map{c} (S^+)\dual\}^{\perp}.\\
\end{array}
$$
Note that $A_U = \bigcap_{v\in U}A_v \subset Q^-$ unless $U$ is 
an $\alpha$-plane $A_s$, in which case $A_U = B_s$---with similar 
remarks holding for $B_U$. If $\dim U = 1,2,3$ then $\dim A_U = \dim B_U = 
4,2,1$ respectively; if $U = B_t$ is a $\beta$-plane then $A_{U}$ is 
1-dimensional and is spanned by $t \in S^-$.  

By triality one can now 
extend the notation of diagram (\ref{quad-pic}) to arbitrary {\it isotropic 
subspaces} $v,s,t$. The following properties are essentially tautological:

\begin{lemm}
\label{AB1}
For any isotropic subspace $U\subset V$, $S^+$ or $S^-$, with $\dim U \not= 3$,
we have:

\itemitem{(i)}
$A_{B_U} = B_{A_U} =U$;
\itemitem{(ii)}
$A_{A_U} = B_U$ and $B_{B_U} = A_U$.
\end{lemm}

\begin{lemm}
\label{AB2}
Suppose that $U,R \subset V\oplus S^+ \oplus S^-$ are 2-dimensional 
isotropic subspaces contained in different summands, and that 
Clifford multiplication vanishes on $U\otimes R$. Then either $U=A_R$, $R=B_U$ 
or $U=B_R$, $R=A_U$.
\end{lemm}

Finally, fix a triality action $\rho :S_3 \rightarrow O(V\oplus S^+ 
\oplus S^-)$
as in \S\ref{group}.

\begin{lemm}
\label{AB3}
For any $g\in S_3$ and isotropic subspace $U\subset V$, $S^+$ or $S^-$ 
we have:
$$
g A_U = \cases{A_{gU}&\cr B_{gU}&\cr}
\qquad
g B_U = \cases{B_{gU}& if ${\rm sgn}\ g = +1$,\cr A_{gU}& if ${\rm sgn}\ 
g = -1$.\cr}
$$
\end{lemm}

\prf
Suppose that $U\subset V$ (the argument for $U\subset S^\pm$ is the same).
Then $A_U = \{ t\in S^- | c(u,s,t)=0\forall s\in S^+, u\in U  \}$
using the trilinear form (\ref{trilinear}). Since this form is $S_3$-invariant
$gA_U$ consists of $x=gt \in V,S^+$ or $S^-$ (depending on $g$) such that 
$c(gu,gs,x) =0$ for all $u\in U$, $s\in S^+$ (where we temporarily 
disregard the order of 
the arguments of $c$), i.e. $c(u',s',x)=0$ for all $u' \in gU$ and $s' \in 
gS^+$. This shows that $gA_U = A_{gU}$ or $B_{gU}$; the assertion 
that which one depends on
${\rm sgn}\ g$ follows easily from diagram (\ref{quad-pic}). 
\qqed

\subsection{Triality action on the Lie algebra}
\label{liealg}

We need next to understand how the Lie algebra 
$\gothg = \so(8)$
decomposes under the action of $S_3$. We write $\gothg = \gothh \oplus
\gothg_+ \oplus \gothg_-$ where $\gothh$ is a Cartan subalgebra
with orthonormal basis $\bolde_1,\bolde_2,\bolde_3,\bolde_4$
and $\gothg_{\pm}$ are spanned by the positive and negative root 
spaces
with respect to simple roots:
$$
\dynk\dddiv{\boldf_1=\bolde_1 - \bolde_2}.{\boldf_0=\bolde_2 
- \bolde_3}.{\boldf_3=\bolde_3 + 
\bolde_4}.
{\boldf_2=\bolde_3 - \bolde_4}.
$$
The 24 roots of $\gothg$ are $\pm 
\bolde_{i} \pm \bolde_{j}$.
The group $S_3$ acts on the Dynkin diagram, and its action on $\gothg$ is 
obtained by 
examining the action on the root spaces. This can be made transparent by 
rewriting the 12 positive roots as:
$$
\begin{array}{llll}
	\boldf_1, & \boldf_0+\boldf_1, & \boldf_0+\boldf_2+\boldf_3, 
	& \boldf_0,  \\
	\boldf_2,  & \boldf_0+\boldf_2, & \boldf_0+\boldf_1+\boldf_3, 
	&\boldf_0+\boldf_1+\boldf_2+\boldf_3,  \\
	\boldf_3,  & \boldf_0+\boldf_3, &\boldf_0+\boldf_1+\boldf_2, & 
	2\boldf_0+\boldf_1+\boldf_2+\boldf_3.
\end{array}
$$

\medskip\noindent
\underline{\it $\sigma$-action}. 
From the Dynkin diagram,  
the involution $\sigma$ acts on $\gothh$ by $\bolde_{4} 
\leftrightarrow -\bolde_{4}$, fixing the other $\bolde_{i}$. Thus
$\gothh = \three^{+} \oplus \one^{-}$ where the superscript denotes 
the eigenvalue of each summand. Alternatively, $\sigma$ interchanges 
$\boldf_2$, $\boldf_3$ and fixes $\boldf_0$, $\boldf_1$. 
From this it follows that $\gothg_{\pm} = 
\six^{+} \oplus (\three^{+} \oplus \three^-)$, and hence
$$
\so(8) = \two\one^+ \oplus \seven^-
\quad
\hbox{under the action of $\sigma$.}
$$

\begin{rem}\rm
Alternatively, one can see this from the fact that $\so(8) \cong\bigwedge^2
V$ where $\sigma$ acts on $V = \c^8$ by reflection in a hyperplane 
$H\subset V$. Then $\two\one^+ = \bigwedge^2 H$ and $\seven^- = H$.
\end{rem}

\medskip\noindent
\underline{\it $\tau$-action}. 
$\tau$ cyclically permutes
$\boldf_1,\boldf_2,\boldf_3$ and fixes
$\boldf_0$, so that 
$\gothh = \two \oplus 
\one^{\omega}\oplus \one^{\omega^2}$ where $\omega = e^{2\pi i/3}$.
For $\gothg_{\pm}$, we observe that 
the entries of each of the first three columns of the array above
are cyclically permuted by $\tau$, while those of the last column are 
invariant. It follows that each of 
$\gothg_{\pm} = \six \oplus \three^{\omega} \oplus \three^{\omega^2}$,
and hence:
$$
\so(8) = \one \four \oplus \seven^{\omega} \oplus \seven^{\omega^2}
\quad
\hbox{under the action of $\tau$.}
$$

\begin{rem}\rm
The invariant 14-dimensional subspace is the Lie algebra of $G_2 = 
\aut \o$. This is well-known, and follows from (\ref{spin8subgp2}).
\end{rem}

\section{Galois $\spin(8)$ bundles}
\label{trigonal}

After some some preliminary remarks on spin bundles in 
\S\ref{spinbundles} we shall introduce the main objects of this paper, 
Galois $\spin(8)$-bundles, in \S\ref{bunS3} (see Definition 
\ref{definition}).

\subsection{Spin bundles}
\label{spinbundles}

We next consider principal 
$G$-bundles
$F\rightarrow X$ where $X$ is a curve and $G= \spin(2n)$ an even 
complex spin group. Given a representation $\rho : G \rightarrow 
SL(W)$ we can form a vector bundle $W_F = F \times_{\rho} W$. 

\begin{lemm}
\label{spin-sum}
If $F_1$ is a $\spin(2m)$-bundle and $F_2$ is a $\spin(2n)$-bundle, 
then 
there is a $\spin(2m+2n)$-bundle $F_1 + F_2$ with
$$
\begin{array}{rcl}
V_{F_1 + F_2} &=& V_{F_1} \oplus V_{F_2},\\
S^+_{F_1 + F_2} &=& S^+_{F_1}\otimes S^+_{F_2} \oplus 
S^-_{F_1}\otimes S^-_{F_2},\\
S^-_{F_1 + F_2} &=& S^+_{F_1}\otimes S^-_{F_2} \oplus 
S^-_{F_1}\otimes S^+_{F_2}.\\
\end{array}
$$
\end{lemm}

\begin{rem}\rm
More generally,
$\nm(F_1 + F_2) = \nm(F_1) \otimes \nm(F_2)$ if $F_1$ and $F_2$ are Clifford 
bundles and $\nm$ is the spinor norm.
\end{rem}

\prf
Since the spin groups are simply connected the inclusion $SO(2m)\times 
SO(2n) \hookrightarrow SO(2m+2n)$, $(g,h) \mapsto g\oplus h$ has a 
unique lift $\lambda$ making the following diagram commute:
$$
\begin{array}{ccc}
	\spin(2m)\times \spin(2n) & \map{\lambda} & \spin(2m+2n)  \\
&&\\
	4:1\downarrow &  & \downarrow 2:1  \\
&&\\
	 SO(2m)\times SO(2n)  & \hookrightarrow & SO(2m+2n).\\
\end{array} 
$$
The bundle $F_1 +F_2$ is obtained by applying $\lambda$ to transition 
functions of $F_1$ and $F_2$ with respect to a sufficiently fine open 
cover of the curve. The identification $V_{F_1 + F_2} =V_{F_1} \oplus 
V_{F_2}$ is then immediate from the construction.

Denote by $V_{1}= N_{1}\oplus N_{1}\dual$ and $V_{2}= N_{2}\oplus 
N_{2}\dual$ 
the orthogonal representations of $\spin(2m)$ and $\spin(2n)$ 
respectively, decomposed into 
maximal isotropic subspaces. Then $N= N_{1}\oplus N_{2}$ is maximal 
isotropic in $V = V_{1}\oplus V_{2}$ and
$$
\begin{array}{ccc}
	\bigwedge^{\rm even}N & = & 
	\bigwedge^{\rm even}N_{1}\otimes \bigwedge^{\rm even}N_{2}
	\oplus 
	 \bigwedge^{\rm odd}N_{1}\otimes \bigwedge^{\rm odd}N_{2}, \\
	\bigwedge^{\rm odd}N & = & \bigwedge^{\rm even}N_{1}\otimes 
	\bigwedge^{\rm odd}N_{2}
	\oplus 
	 \bigwedge^{\rm odd}N_{1}\otimes \bigwedge^{\rm even}N_{2}.
\end{array}
$$
From this it follows that the spinor bundles of $F_1 + F_2$ are as 
asserted.
\qqed

If $P\subset G$ is a parabolic subgroup we can form the bundle of 
homogeneous 
spaces $F/P$ with fibre $G/P$. The $G$-bundle $F$ is said to be {\it 
stable} (resp.{\it semistable}) if for every maximal parabolic 
$P\subset G$ and every section $s: X \rightarrow F/P$ one has 
$$
\deg s^* T^{\rm vert}_{F/P} > 0 
\qquad
\hbox{(resp. $\geq 0$)}
$$
where $T^{\rm vert}$ denotes the vertical tangent bundle. (See 
\cite{Rth}.)
In the case of $G= \spin(2n)$ there are (up to conjugacy) $n$ 
maximal parabolics and the spaces $G/P$ are the grassmannians of 
isotropic subspaces 
$U\subset V$ of dimension $d=1,\ldots, n-2$ and the two spinor 
varieties of isotropic subspaces of dimension $d=n$. A section $s: X 
\rightarrow F/P$ 
corresponds to an isotropic subbundle $\underline U \subset V_F$ and 
it is easy to check that the stability condition reduces to the slope 
inequality 
$$
\mu(\underline U) < \mu(V_F)
\qquad
\hbox{
$\forall$ isotropic subbundles $\underline U \subset V_F$.}
$$
In other words, $F$ is a (semi)stable $\spin(2n)$-bundle if and only 
if $V_F$ 
is (semi)stable as an orthogonal vector bundle. (See also \cite{O}, 
Lemma 1.2.) 

Recall that, just as for vector bundles, the moduli problem for spin 
bundles requires an equivalence relation on semistable bundles 
coarser than isomorphism, called S-equivalence. For the general 
definition in the context of principal bundles we refer to \cite{KNR};
for our purposes the following remarks will be sufficient.  

Suppose that $W$ is an orthogonal vector bundle with a 
filtration
$$
U_1 \subset U_2 \subset \cdots \subset U_k \subset U_k^{\perp}\subset \cdots
U_2^{\perp}\subset U_1^{\perp}\subset W,
$$
where $U_{1},\ldots,U_k$ are destabilising subbundles (i.e. isotropic 
of degree 0) and each `quotient' 
$$
(U_{i+1}/U_i)\oplus (U_{i}^{\perp}/U_{i+1}^{\perp})
\cong 
(U_{i+1}/U_i)\oplus (U_{i+1}/U_{i})\dual
$$
is a stable orthogonal bundle, with its natural orthogonal structure 
(we include $U_0 = 0$). Clearly such a filtration always exists, and 
$k=0$ if and only if $W$ is stable. The graded orthogonal bundle
$$
\gr(W) =  \bigoplus_{i=0}^{k-1}\left\{ (U_{i+1}/U_i)\oplus (U_{i+1}/U_{i})\dual
\right\}\oplus U_k^{\perp}/U_k
$$
is independent of the filtration, and two orthogonal bundles $W,W'$ of the 
same rank are said to be S-equivalent if and only if their 
graded bundles are isomorphic as orthogonal bundles:
$$
W \sim W' 
\qquad \Longleftrightarrow  \qquad
\gr(W) \cong \gr(W').
$$

Recall that by \cite{R} Proposition 4.5 $V_F$ is stable as an 
orthogonal bundle if and only if $V_F$ is polystable as a 
vector bundle, i.e.
$$
V_F = V_1 \oplus \cdots \oplus V_k,
$$
where the summands $V_i$ are stable as vector bundles and nonisomorphic.

Let us now restrict to the case of $\spin (8)$,
and suppose that the spinor bundles are polystable of the same shape, 
that is, with stable summands of the same ranks as the $V_i$:
$$
\begin{array}{rcl}
S^+_F &=& S^+_1 \oplus \cdots \oplus S^+_k,\\
S^-_F &=& S^-_1 \oplus \cdots \oplus S^-_k.\\
\end{array}
$$
In this situation we shall need later on to understand the automorphism group 
of the spin bundle $F$. This is determined by its action on the three 
bundles above; by stability and orthogonality the automorphism group of
each of these is $\mu_2 ^k$ where $\mu_2 = \{\pm 1\}$. $\aut F$ is then 
described by the following diagram with exact rows:
\begin{equation}
\label{aut}
\begin{array}{ccccccc}
0\rightarrow &\mu_2& \map{\beta} & \aut F & \map{\alpha} & (\mu_2)^k &
\rightarrow 0\\
&&&&&&\\
&\|&&\qquad \uparrow {\rm inclusion}&&\qquad \uparrow {\rm diagonal}&\\
&&&&&&\\
0\rightarrow &\mu_2& \rightarrow & Z(\spin(8)) & \rightarrow & \mu_2 &
\rightarrow 0.\\
\end{array}
\end{equation}
The map $\alpha$ is the representation on $V_{F}$, while $\beta(-1)$ 
acts as $-1$ on each of $S_F^{\pm}$ and as $+1$ on $V_F$.
Explicitly, the elements of $\aut F$ are the $4\times 2^{k-1} = 2^{k+1}$
matrices (acting on 
$V_F\oplus S_F^+ \oplus S_F^-$)
\begin{equation}
\label{split-autos}
\pmatrix{\varepsilon&&\cr &\varepsilon& \cr &&\varepsilon\cr},
\pmatrix{\varepsilon&&\cr &-\varepsilon& \cr &&-\varepsilon\cr},
\pmatrix{-\varepsilon&&\cr &\varepsilon& \cr &&-\varepsilon\cr},
\pmatrix{-\varepsilon&&\cr &-\varepsilon& \cr &&\varepsilon\cr},
\end{equation}
where $\varepsilon \in (\mu_2)^k$ satisfies $\prod \varepsilon_{i} 
=1$, i.e. contains an even number of $-1$s (and in each matrix 
denotes a diagonal block).

\subsection{Bundles with $S_3$-action}
\label{bunS3}

Let $X$ be a smooth curve acted on by the group $S_3$. We shall assume 
that the action is faithful and
that each involution has nonempty fixed-point set.

The group $S_3$ acts in 
two ways on isomorphism classes of principal 
$\spin(8)$-bundles on $X$. First, it has a right 
`triality' action $F\mapsto F^{u}$ (where $u \in S_3$) by 
outer 
automorphisms of the structure group. These are defined up to inner 
automorphisms, as in Lemma \ref{outerautos}, 
but inner automorphisms preserve the isomorphism class 
of the bundle. Second, the right action of $S_3$ on $X$ induces 
by pull-back a left action on bundles, $F\mapsto u(F) := (u^{-1})^* F $. 
We shall be interested in (semistable) 
bundles for which these two actions agree, i.e. 
the fixed-point set $\ff_X \subset \mm_X(\spin(8))$ of the 
$S_3$-action $F\mapsto u^* F^u$, $u\in S_3$.

If $F\in \ff_X$ then there is a nontrivial {\it Mumford group} $\gg_F$ 
consisting of pairs $g\in S_3,\ \lambda:g(F)\idfy F^g$ with 
multiplication law
$$
(g,\lambda)(h,\mu) = (gh, \lambda^h \circ g\mu).
$$ 
Here $\lambda^h : g(F^h)=g(F)^h \idfy F^{gh}$ and $g\mu : g h(F) \idfy 
g(F^h)$ are the natural induced isomorphisms.
This group is an extension
\begin{equation}
\label{mumford}
1\longrightarrow \aut F
\longrightarrow \gg_F
\longrightarrow S_3
\longrightarrow 1.
\end{equation}
Note in particular that there is an $S_3$-action on $\aut F$ induced, 
via conjugation in $\gg_F$, by this sequence. It is easy to verify 
that the invariant subgroup $(\aut F)^{S_3}$ consists of $\alpha \in \aut 
F$ commuting with $S_3$ in the sense that for all $g\in S_3 $ and 
$\lambda: g(F)\idfy F^g$ the 
following diagram commutes:
$$
\begin{array}{ccc}
	g F & \map{\lambda} & F^g  \\
	&&\\
	g\alpha \downarrow &  & \downarrow \alpha^g  \\
	&&\\
	g F & \map{\lambda} & F^g.\\
\end{array}
$$

\begin{rem}\rm
\label{S3-auto-remark}
If $F$ is stable as a
$\spin(8)$ bundle then the vector bundles $V_F, S_F^{\pm}$ are polystable
with $k$ stable summands, say. The automorphism group $\aut F$ is then 
the group of order $2^{k+1}$ described in the last section 
((\ref{aut}) and (\ref{split-autos})).
The action of $S_3$ on $\aut F$ determined by the sequence 
(\ref{mumford}) permutes, for each $\varepsilon \in (\mu_2)^k$, 
$\prod \varepsilon_{i} 
=1$, the last three matrices of 
(\ref{split-autos})) in the natural way. In particular we have
\begin{equation}
\label{S3-autos}
(\aut F )^{S_3} = \ker \{\Pi : (\mu_2)^k \rightarrow \mu_2 \}.
\end{equation}
\end{rem}

\begin{lemm} 
\label{galoisbundle}
Suppose $F\in \ff_{X}$. Then the following data are equivalent:
\begin{enumerate}
\item
A splitting $\Lambda : S_3 \rightarrow \gg_F$ of the Mumford 
sequence (\ref{mumford}).
\item
A lift of the $S_3$-action on $X$ to the principal bundle $F$, in the 
sense that there is for each $u\in S_3$ a commutative diagram
$$
\begin{array}{ccc}
	F & \map{\lambda(u)} & F^g  \\
	&&\\
	\downarrow &  & \downarrow  \\
	&&\\
	X & \map{g\in S_3} & X
\end{array}
$$
and that $\lambda(u)^w \circ \lambda(w) = \lambda(uw)$
for all $u,w \in S_3$.
\item
An orthogonal lift of the $S_3$-action on $X$
to the {\it Chevalley bundle} $V_F \oplus S^+_F \oplus S^-_F$ such 
that:

\itemitem{(i)}
the 3-cycle $\tau\in S_3$ permutes the summands cyclically $V_F \idfy 
\tau^*S^+_F \idfy (\tau^2)^*S^-_F \idfy V_F$;
\itemitem{(ii)} 
the involution $\sigma \in S_3$ lifts to an involution $V_F \idfy 
\sigma^* V_F$,
acting in the fibre $(V_F)_x$ at a fixed point $x\in X$ of $\sigma$ 
as reflection in a hyperplane on which the quadratic form is nondegenerate, 
and exchanges the spinor bundles $S^+_F \idfy \sigma^*S^-_F$.
\end{enumerate}
\end{lemm}

\prf
The equivalence of 1 and 2 is essentially obvious and we leave
its verification to the reader.
To show that 2 is equivalent to 3
we shall represent the bundle $F$ by a $\spin(8)$-valued 
Cech cocycle $\{g_{ij}\}$ with respect 
to an open cover $\{U_i\}$ of the curve. We can assume that this open 
cover 
is invariant under the action of the finite group $S_3$, and we shall 
denote the image of $U_i$ under $w\in S_3$ by $U_{i^w}$. 

If $\{h_{ij}\}$ represents a second 
bundle $F'$ then a bundle isomorphism $F\idfy F'$ is represented by a 
cochain 
$\{f_i\}$ satisfying
$$
f_i g_{ij} = h_{ij} f_j.
$$
For each $w\in S_3$ the bundle $w^* F$ has cocycle $\{w^*(g_{i^w 
j^w})\}$; 
while the bundle $F^{w^{-1}}$ has cocycle $\{w g_{ij} w^{-1}\}$, where the 
conjugation
takes place in the semidirect product $\spin(8) \sdp S_3 \subset 
O(V\oplus S^+ \oplus S^-)$. The condition that $F^{w^{-1}} \cong w^* 
F$ is therefore:

\begin{equation}
\label{gal-condition1}
\hbox{\it $\exists$ $\spin(8)$-valued cochain 
$\{f_i\}$ satisfying
$f_i^{-1} w^*(g_{i^w j^w}) f_j = w g_{ij} w^{-1}$.
}
\end{equation}

Next now suppose that the vector bundle $V_F \oplus 
S^+_F \oplus 
S^-_F$ admits an orthogonal lift of the $S_3$-action. It is easy to 
check 
that necessary and sufficient conditions for this are 
the existence, for each $w\in S_3$, 
of an $O(V \oplus S^+ \oplus S^-)$-valued
cochain $\{\ftil_i \}$ satisfying

\begin{equation}
\label{gal-condition2}
\ftil_i g_{ij} = w^*(g_{i^w j^w}) \ftil_j,
\end{equation}
together with suitable compatibility 
assumptions on the cochains $\{\ftil_i \}$ as $w\in S_3$ varies. 
We shall show that (\ref{gal-condition2}) reduces to (\ref{gal-condition1}) 
when the lift
satisfies the conditions (i) and (ii) of the lemma. 
Namely, these conditions
are equivalent to requiring that $\{\ftil_i \}$ takes values in 
$\spin(8) \sdp S_3 \subset O(V\oplus S^+ \oplus S^-)$ and projects to 
the 
element $w^{-1}\in S_3$. (Note that the 1-dimensional $-1$-eigenspaces 
at the fixed points when $w$ is an involution 
(3 part (ii) of the lemma) arise from the choices
of reflections of $V\oplus S^+ \oplus S^-$ needed to define $\rho: S_3
\hookrightarrow O(V\oplus S^+ \oplus S^-)$.)

Consequently, for each $w\in S_3$ the cochain defining the lift of 
$w^{-1}$ 
to the Chevalley bundle can be expressed as $\ftil_i = (f_i,w^{-1})$
where $f_i$ is $\spin(8)$-valued. Then (using (\ref{sdp}) in 
\S\ref{spin})
$$
\begin{array}{rcl}
\ftil_i g_{ij}& =& (f_i,w^{-1})(g_{ij},1)\\
&=& (f_i w g_{ij}w^{-1},w^{-1})\\
\end{array}
$$
while
$$
\begin{array}{rcl}
w^*(g_{i^w j^w}) \ftil_j & =&(w^*(g_{i^w j^w}),1)(f_j,w^{-1}) \\
&=& (w^*(g_{i^w j^w})f_j,w^{-1}).\\
\end{array}
$$
So we see that (\ref{gal-condition2}) reduces to 
(\ref{gal-condition1}) as required.
\qqed

\begin{define}
\label{definition}
\itemitem{(i)}
By a {\it Galois $\spin(8)$-bundle} we shall mean a pair $(F,\Lambda)$ 
where $F$ is a principal $\spin(8)$-bundle on a curve $X$ on which 
$S_3$ acts faithfully, and $\Lambda: S_3 \rightarrow \gg_F$ is a 
lift of the group action to $F$ in the sense of Lemma \ref{galoisbundle}. 
Semistability and S-equivalence of Galois 
bundles will refer to the corresponding properties of the underlying 
bundles $F$.
\itemitem{(ii)}
Let $\nn_X$ be the set of S-equivalence classes of semistable Galois 
$\spin(8)$-bundles; and let $\ff_X \subset \mm_X(\spin(8))$ denote the 
fixed point set of the 
$S_3$-action $F\mapsto u^* F^u$, $u\in S_3$. For $F\in \ff_X$ we 
denote by $\nn(F)$ the fibre of the forgetful map
$$
\nn_X \rightarrow \ff_X\subset \mm_X(\spin(8)).
$$
\end{define}

The set $\nn(F)$ of Galois structures, if nonempty, is a torsor over the 
cohomology group $H^1(S_3, \aut F)$.

\begin{lemm} {\capit (S.M.J. Wilson)}
\label{wilson}
Let $G$ be a finite group and $H \subset G$ be a normal subgroup of 
coprime order and index. Suppose that $G$ acts on an abelian group
$A$ such that $H,A$ have coprime order. Then for all $q\geq 1$:
$$
\begin{array}{rcl}
	H^{q}(H,A) & = & 0,  \\
	H^{q}(G,A) & \cong & H^{q}(G/H,A^{H}).\\
\end{array}
$$
\end{lemm}

\prf
As a $\z$-module $H^{q}(H,A)$ is annihilated by $|A|$ trivially, and 
is also annihilated by $|H|$ since every element has order dividing 
$|H|$ by \cite{Mac} p.117 Proposition 5.3. $H^{q}(H,A)$ is therefore 
annihilated by ${\rm gcd}(|H|,|A|) = 1$, proving the first part. It 
then follows from this and the Hochschild-Serre spectral sequence that 
for all $q\geq 1$ the inflation maps $H^{q}(G/H,A^{H}) \rightarrow
H^{q}(G,A)$ are isomorphisms (see \cite{Mac} p.355 Exercise 3).
\qqed

In our situation $G=S_3$, $H = \<\tau\>$, $A = \aut F$, the 
quotient $G/H$ acts trivially on $A^H$, and it follows that 
$$
H^q(S_3, \aut F) \cong \ker \{\Pi : (\mu_2)^k \rightarrow \mu_2 \}
\qquad
\hbox{for $q=1,2$.}
$$
From this we conclude:

\begin{prop}
\label{howmanylifts}
Suppose $F \in \ff_X$ is stable with $k$ summands (as a polystable 
vector bundle). 
\itemitem{(i)}
If $k=1$ then there exists a unique Galois structure on $F$.
\itemitem{(ii)} 
If $k>1$ then $H^2(S_3,\aut F) \cong (\mu_2)^{k-1}$. If the Mumford 
sequence (\ref{mumford}) splits then the set $\nn(F)$ of Galois 
structures on $F$ has cardinality $2^{k-1}$.
\end{prop}

\begin{rem}\rm
If $V_F$, and hence also $S_F^{\pm}$, are stable as 
vector bundles (case (i) above) then one can see the uniqueness of the 
Galois structure directly. Consider isomorphisms $\alpha, \beta:V_F \idfy 
\sigma^* V_F$. By stability these differ by a scalar $\alpha = 
\lambda \beta$, $\lambda \in 
\c$, but at a fixed point both $\alpha, \beta$ act with 1-dimensional 
$-1$-eigenspace in the fibre, and this forces $\lambda = 1$.
\end{rem}

\subsection{Further remarks}
\label{hypcase}

We should ask at this point how one can construct Galois bundles. 
Note that since $F$ is determined by its 
spinor bundles $S_F^{\pm}$ it is determined by the orthogonal bundle 
$V_F$:
$$
S^+_F = \tau(V_F),
\qquad
S^-_F = \tau^2(V_F).
$$
(Because of the 
relation $\sigma \tau = \tau^2 \sigma$ (recall \ref{notation} for notation) 
these two bundles are necessarily
interchanged by $\sigma$.)
In the following lemma, which plays a key role in later sections, we
construct a `split' orthogonal bundle $V_F$ which gives rise to a Galois 
spin bundle. We make use of the quotient 
$$
\pi: X \map{2:1} Y = X/\sigma.
$$

\begin{lemm-def}
\label{su2trigonal}
Let $E\rightarrow Y = X/\sigma$ be a rank 2 vector bundle with $\det 
E = 
\oo_Y$.
Then there is a unique Galois $\spin(8)$-bundle $F=F_E\rightarrow 
X$ with
$$
V_{F} = \c^2 \otimes \pi^{*}E \oplus \tau(\pi^{*}E)\otimes 
\tau^2(\pi^{*}E).
$$
\end{lemm-def}

\prf 
Let $F_1 = \c^2 \otimes E'$ and $F_2 = \tau(E')\otimes 
\tau^2(E')$ for any rank 2 bundle $E'\rightarrow X$ with trivial 
determinant.
Each of $F_1,F_2$ is a $\spin(4)$-bundle whose spinor bundles are its 
two 
factors.
So by Lemma \ref{spin-sum}, $F_1 \oplus F_2$ is the orthogonal 
representation of a $\spin(8)$ bundle $F$, and we have:
$$
\begin{array}{rcl}
V_{F} &=& F_1 \oplus F_2,\\
S^+_{F} &=&  \c^2 \otimes \tau(E')\oplus \tau^2(E')\otimes 
E' = \tau(F_1 \oplus F_2),\\
S^-_{F} &=& \c^2 \otimes \tau^2(E') \oplus E' \otimes 
\tau(E') = \tau^2(F_1 \oplus F_2).\\
\end{array}
$$
It follows that $\sigma$ lifts to an isomorphism 
$S_{F}^{+}\idfy \sigma^{*} S_{F}^{-}$. 
Moreover, $\sigma$ lifts to an isomorphism $F_2 \idfy \sigma^* F_2$ 
having in the fibre at each fixed point $x\in X$ one-dimensional  
$-1$-eigenspace $\bigwedge^2 (\tau E')_x$. This extends to an 
isomorphism $V_F \idfy \sigma^* V_F$ with the same property provided
$E'$ is the pull-back of some bundle $E\rightarrow Y$; so we are done.
Note that in this construction we must have $\det E = \oo_Y$ since 
$\det E' = \oo_X$:
by hypothesis $\sigma$ has nonempty fixed-point set,  
so $\pi^*$ is 
injective on line bundles.
\qqed

\begin{rem}\rm
There are, of course, two more quotients of $X$ by choosing 
$\sigma\tau$ or $\sigma\tau^2$ instead of $\sigma$. In fact these are 
both isomorphic to $Y$: the $S_3$ action on $X$ induces an inclusion 
$X \hookrightarrow Y\times Y\times Y$, and the three quotients are 
then simply the projections on the three factors. Under the above 
construction, choosing a different projection has the effect of 
permuting 
the bundles $V_F,S_F^{\pm}$.
\end{rem}

The existence of Galois bundles not of the form \ref{su2trigonal}
(and stable, in fact---note that the bundle
given by \ref{su2trigonal}) is semistable but not stable:
$\pi^* E \hookrightarrow V_F$ is a 
destabilising subbundle) will be seen in section \ref{recillas}.

\medskip

Our next remark is a numerical observation that will reappear in later
(see Remark \ref{dimcount}(ii)).

We might relax the condition that $S_3$ acts faithfully on the curve $X$.
Suppose, in particular, that the normal subgroup $\< \tau \> \subset S_3$ acts 
trivially, i.e. that $S_3$ acts by a single involution of $X$. Then  
a Galois $\spin(8)$-bundle $F$ requires isomorphisms $V_F \idfy 
S^{\pm}_F$; let us interpret this as asking simply for a rank 8 orthogonal 
bundle $V_F \rightarrow X$ together with a lift of the involution 
$\sigma : X\leftrightarrow X$ which at fixed points acts in the fibre 
with 1-dimensional $-1$-eigenspace.

Suppose in particular that $X\rightarrow \bp^1$ is a hyperelliptic 
curve and $\sigma$ the sheet-involution. Then provided $X$ has genus 
$g\geq 4$ there is a 
projective moduli space $\mm$ of such $SO(8)$-bundles which is 
described 
explicitly in \cite{R} Theorem 1. 
This says that 
$$
\mm \cong \grass_{g-4}(Q_1 \cap Q_2) /(\z/2)^{2g+2},
$$
where $Q_1,Q_2 \subset \bp^{2g+1}$ are quadrics, and $\grass_{g-4}$ 
denotes the grassmannian of projective $g-4$-planes isotropic for 
both quadrics. 
In particular the dimension is easy to compute: $\mm$ is cut out in 
the grassmannian by a pair of sections of the bundle $S^2 
U\dual$ where $U$ is the tautological bundle. Hence
$$
\begin{array}{rcl}
	\dim \mm & = & (g-3)(g+5)-(g-3)(g-2)  \\
	 & = & 7(g-3).
\end{array}
$$

\medskip

Finally, we make some remarks concerning 
curves with $\z/3$-action---in particular, for example, cyclic trigonal
curves. In this case we can certainly imitate the above constructions replacing $S_3$ by $\z/3$, with triality action on $\spin(8)$ determined by making a 
choice (there are two) of embedding $\z/3 \hookrightarrow S_3$. The conjugation
action of $S_3$ 
on $\mm_X(\spin(8))$ then restricts to $\z/3$ and we can consider
the fixed-point-set $\ff_X^{\z/3} \subset \mm_X(\spin(8))$. A lifting
of the $\z/3$-action is determined by a splitting of the `restricted'
Mumford group 
$$
\ses{\aut F}{\gg_F^{\z/3}}{\z/3}.
$$
However, in this case the sequence {\it always splits uniquely} since by Lemma
\ref{wilson} $H^q(\z/3,\aut F) =0$ for $q=1,2$. So for our moduli space of
$\z/3$-Galois bundles we can simply take $\nn_X = \ff_X$. Geometrically
a Galois bundle is now a $\spin(8)$-bundle $F$ together with a lift of the 
(chosen triality) $\z/3$-action in the sense of \ref{galoisbundle}(2), 
or equivalently \ref{galoisbundle}(3) omitting the requirement (ii).

We expect that the results outlined in the introduction for the $S_3$ case
should hold also for $\z/3$. Indeed, the discussion of \S\ref{local} 
in the next 
section goes through unimpeded, as does the dimension calculation of the 
moduli space in \S\ref{gen-dim} (see Remark \ref{cyclic-dim}). However, 
the difficulty arises in computing the semistable boundary (\S\ref{boundary}),
where we make essential use of the elements of order 2 in $S_3$. Possibly
one could get round this and prove a corresponding result for $\z/3$, 
but we have not pursued the question here.

\section{The moduli space}
\label{moduli}

We wish to construct a moduli space $\nn_X$ of Galois $\spin(8)$ 
bundles on our $S_3$-curve $X$. 
Using Proposition \ref{howmanylifts} this can be modelled on
the fixed-point set $\ff_{X} \subset \mm_{X}(\spin(8))$
of the `conjugation' action 
$F\mapsto u^{-1}F^{u}=  u^* F^{u}$, $u\in S_3$.
We shall show in this section that $ \nn_X$ inherits from $\ff_X$ 
the structure of an analytic space which is smooth over stable
spin bundles. $\nn_X$ and $\ff_X$ are locally isomorphic at stable 
vector bundles (i.e. at $F$ such that $V_F$ is stable as a vector 
bundle), while at polystable vector bundles $\nn_X$ resolves normal 
crossing singularities in $\ff_X$:
 \begin{center}
\leavevmode
\epsfxsize=3.2in\epsfbox{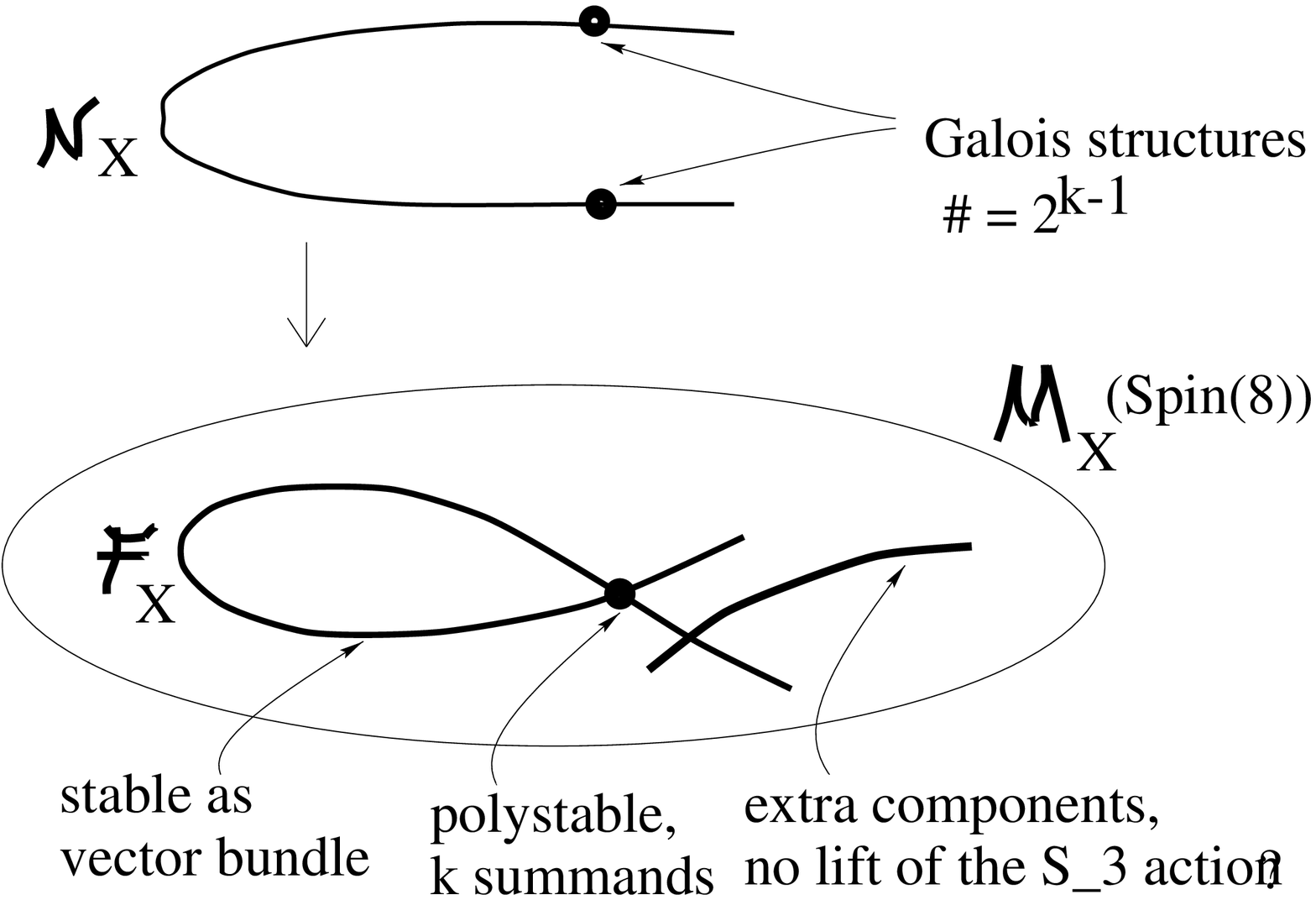}
\end{center}
We analyse the semistable boundary (i.e. Galois bundles nonstable as 
spin bundles) and show that, with the additional assumption that
$X/S_3 \cong \bp^1$, this 
consists precisely of equivalence classes of bundles of the form 
\ref{su2trigonal}.

\subsection{Local moduli}
\label{local}

To begin, we examine the $S_3$-action on the Kodaira-Spencer map at a 
stable spin bundle in $\mm_X(\spin(8))$.

Suppose that $F$ is a principal $G$-bundle $F$, for some reductive group $G$,
represented by a Cech cocyle (transition functions) $\{g_{ij}\}$ with 
respect an open cover $\{U_i \}$ of $X$. If $\ad F$ is the vector 
bundle $F\times_{\ad} \gothg$ then a cohomology class $\xi \in 
H^1(X,\ad F)$ is represented by a cocycle $\{ \xi_{ij} \}$ of 
$\gothg$-valued functions on $U_i \cap U_j$ satisfying
$$
	0  = \xi_{ij}+\xi_{ji}^{g_{ij}},   \quad
	 0 = \xi_{ij}+\xi_{jk}^{g_{kj}}+ \xi_{ki}^{g_{ik}g_{kj}},
\qquad
\forall\ i,j,k;
$$
where $\xi^g := {\rm Ad}(g)\xi$ for $\xi\in \gothg$ and $g\in G$. 
Exponentiating these conditions, i.e. applying $\exp: \gothg 
\rightarrow G$, one gets precisely the cocycle conditions for 
transition functions $\{g_{ij}\exp(\xi_{ij}) \}$. Denote the 
corresponding $G$-bundle by $F_{\xi}$. If the bundle $F$ is stable 
then so is $F_{\xi}$ in a neighbourhood of $0\in H^1(X,\ad F)$, and we 
have a rational map
$$
\ks : H^1(X,\ad F) \rightarrow \mm_{X}(G),
$$
mapping $0\mapsto F$ and complete in the sense of \cite{Rth} Theorem 
4.2. That is, $\ks$ induces an isomorphism from a neighbourhood of $0\in H^1(X,\ad F)
/\Gamma_F$, where $\Gamma_F = \aut F/Z(G)$, to a neighbourhood of 
$F\in \mm_X(G)$.

Now let $G = \spin(8)$, and suppose that $(F,\Lambda)$ is a Galois 
bundle (Definition \ref{definition}). Since the adjoint representation 
of $\spin(8)$ is preserved by triality it follows from Lemma 
\ref{galoisbundle}(2) that $\Lambda$ induces an orthogonal lift of the 
$S_3$-action on $X$ to the vector bundle $\ad F$. (This is completely 
analogous to the lift of \ref{galoisbundle}(3) to the Chevalley 
bundle; and indeed can be expressed in terms of \ref{galoisbundle}(3),
since as $\spin(8)$-bundles $\ad F = \bigwedge^{2}V_F \cong 
\bigwedge^{2}S^+_F \cong\bigwedge^{2}S^-_F$.)
In particular $\Lambda$ induces an action of $S_3$ on $H^1(X, \ad F)$.

\begin{prop}
Suppose that $S_3$ acts faithfully on the curve $X$ and that
$(F,\Lambda)$ (where $\Lambda \in \nn(F)$) 
is a stable Galois $\spin(8)$-bundle. Then the Kodaira-Spencer map
$\ks$ is $S_3$-equivariant, where the $S_3$-action on $H^1(X,\ad F)$ 
is that induced by $\Lambda$ and the action on 
$\mm_X(\spin(8))$ is $E\mapsto u^{*}E^u$ for $u\in S_3$. 
\end{prop}

The proof is a straightforward calculation with transition functions, 
using (\ref{gal-condition1}) and the cocycle conditions above for 
elements $\xi \in H^{1}(X,\ad F)$, and we omit the details. 

For each $(F,\Lambda) \in \nn_X$ we shall denote the $S_3$-invariant 
subspace of $H^{1}(X,\ad F)$ by $H^{1}(X,\ad F)^{\Lambda}$. 
It follows from the proposition that the Kodaira-Spencer map restricts 
to a map
$$
\ks_{\Lambda}: H^{1}(X,\ad F)^{\Lambda} \rightarrow \ff_X.
$$
Suppose 
first that $V_F$ is a stable vector bundle, so that (by Proposition 
\ref{howmanylifts})
$\Lambda$ is unique.
The map $\ks_{\Lambda}$ factors to an injective map on the quotient 
$H^{1}(X,\ad F)^{\Lambda}/\Gamma_F^{\Lambda}$ where $\Gamma_F^{\Lambda}
= (\aut F)^{S_3}/Z(\spin(8))^{S_3}$. But from Remark \ref{S3-auto-remark} 
we know 
that $Z(\spin(8))^{S_3}= \{1\}$; and when $V_F$ is a stable vector bundle
$(\aut F)^{S_{3}}$ is also trivial. 
We therefore have a commutative diagram (more precisely, the 
horizontal maps are defined on neighbourhoods of~0):
$$
\begin{array}{ccc}
H^1(X,\ad F)/\Gamma_F & \hookrightarrow & \mm_X(\spin(8)) \\
&&\\
\uparrow && \uparrow \\
&&\\
H^1(X,\ad F)^{\Lambda} &\hookrightarrow & \ff_X. \\
\end{array}
$$
The spaces $\nn_{X}$ and $\ff_X$ are therefore locally isomorphic, 
determining a smooth analytic structure on $\nn_X$ with tangent space
$H^1(X,\ad F)^{\Lambda}$.

\begin{rem}\rm
In \S\ref{recillas} we shall construct (for the case  
$X/S_3 \cong \bp^1$) examples of
$(F,\Lambda) \in \nn_X$ for which the vector bundle $V_F$ is stable.
\end{rem}

More generally, for polystable $V_F$ with $k$ summands, consider the 
set of Galois structures $\nn(F) 
=\{\Lambda_{1},\ldots,\Lambda_{2^{k-1}}\}$. This set is acted on 
transitively and faithfully by the group $(\aut F)^{S_{3}}$. 
The subgroup of $\Gamma_F$ 
preserving each $H^1(X,\ad F)^{\Lambda_i}$ is therefore trivial and 
so we have inclusions (from a neighbourhood of 0 to a 
neighbourhood of $F$):
$$
H^{1}(X,\ad F)^{\Lambda_i} \hookrightarrow \ff_X;
\qquad
i=1,\ldots,2^{k-1}.
$$
Each space can be viewed as Zariski tangent space of $\nn_X$ at 
$(F,\Lambda_i)$, and we conclude:

\begin{prop} 
\label{tgt-space}
The moduli space $\nn_X$ is a complex analytic space mapping 
holomorphically to $\ff_X\subset \mm_X(\spin(8))$.
At a stable Galois bundle $(F,\Lambda)$ it has Zariski tangent space
$$
T_{F,\Lambda} \nn_X  = H^{1}(X,\ad F)^{\Lambda}.
$$
\end{prop}

We shall compute the dimension of the tangent space in section 
\ref{gen-dim}, and deduce (in the case of interest $X/S_3 \cong \bp^1$) 
that $\nn_X$ is 
smooth at stable points (Theorem \ref{dimtheorem}).

\subsection{The semistable boundary}
\label{boundary}

Our next task is to describe Galois bundles which are semistable but 
not stable. We shall assume from now on that the quotient $X/S_3$ 
is isomorphic to $\bp^1$, and under this assumption
it turns out that up to S-equivalence
these bundles are exactly those described in Lemma
\ref{su2trigonal}. To recall, let $E$ be a rank 2 vector bundle on $X$ 
pulled back from $X/\sigma$. 
We consider the Galois $\spin(8)$-bundle $F_E = 
(V_{E},S^+_E,S_E^-)$ where
\begin{equation}
\label{split-form}
\begin{array}{rcl}
	V_E & = & (\tau E)(\tau^2 E)\oplus E \c^2 ,\\
	S^+_E & = & (\tau^2 E)E \oplus (\tau E) \c^2 , \\
	S^-_E & = & E(\tau E) \oplus (\tau^2 E)\c^2 . \\
\end{array}
\end{equation}
We have written $E\c^2 = E\otimes \c^2$ etc for brevity. Note that if a 
trivialisation of the determinant 
line bundle $\det E$ is fixed then fibrewise we are in 
the situation of \S\ref{2by2}: each rank 2 factor is naturally 
self-dual, taking adjugates is a well-defined involution of each rank 
4 summand, and we have an explicit description of the Clifford 
multiplications among the three rank 8 bundles.

\begin{theo}
\label{ssboundary}
Suppose that $S_3$ acts faithfully on a curve $X$
with quotient $X/S_3 \cong \bp^1$. 
A semistable Galois $\spin(8)$ bundle $F\rightarrow X$ is unstable if 
and only if it is of the form $F = F_E$ for some $E\in \ss\uu_{X/\sigma}(2)$ 
up to S-equivalence.
\end{theo}

We have seen in \S\ref{spinbundles} that stability of $F$ is 
equivalent to stability of $V_F$ as an orthogonal bundle; so suppose 
that $U\subset V_F$ is a destabilising subbundle. This means that $U$ 
is an isotropic subbundle of degree 0. In particular it is semistable 
as a vector bundle and $\rk U \leq 4$.

Using the lift of the $S_3$-action on $X$ to the Chevalley bundle 
$W_F = V_F \oplus S_F^+ \oplus S_F^-$ we can consider the image 
bundles $gU \subset W_F$ for group elements $g\in S_3$. In addition, 
we can consider the associated bundles (see \S\ref{6quadrics})
$$
\begin{array}{rcl}
A_U &=& \im \{ U\otimes S_F^+ \rightarrow S_F^-\}^{\perp},\\
B_U &=& \im \{ U\otimes S_F^- \rightarrow S_F^+\}^{\perp};\\
\end{array}
$$
as well as combinations of these constructions such as $A_{A_{U}}$, 
$B_{gU}$ etc. We shall keep track of these subbundles of the 
Chevalley bundle using the results of \S\ref{6quadrics}. 

To begin, suppose that $U\subset V_F$ has  
$\rk U = 4$.
Then one of $A_U\subset S_F^-$ or $B_U \subset S_F^+$ is a 
destabilising line 
subbundle, and so by the action of $\tau \in S_3$ we get a 
destabilising 
line subbundle in $V_F$. 

If, on the other hand, $\rk U =3$ then it is---at least locally---the 
intersection 
of a pair of rank 4 isotropic subbundles parametrising isotropic 
subspaces of opposite families in the fibres. Globally we obtain a 
rank 4 subbundle of the pull-back of $V_F$ to an \'etale double cover 
of $X$; however, the class of this double cover is the Stiefel-Whitney 
class $w_1(F)$, which is zero since $F$ is spin. We therefore have a 
pair of rank 4 isotropic subbundles globally on $X$. 
Since $U$ and $V_F$ are semistable these rank 4 bundles have degree 0 
and so 
again $V_F$ has (two) destabilising line subbundles.

It is therefore sufficient to assume that $\rk U \leq 2$. We shall 
consider for any subbundle $U\subset V_F$ the Clifford 
multiplication
$$
\begin{array}{rcc}
\tau\sigma U \otimes \tau^2 U & \hookrightarrow & S_F^+ \otimes S_F^- \\
&&\\
&\mu \searrow & \downarrow \\
&&\\
&&V_F\\
\end{array}
$$

\begin{lemm}
\label{mu=0}
If $\rk U =1$ then $\mu =0$; if $\rk U =2$ then $\rk \mu \leq 1$.
\end{lemm}

\begin{rem}\rm
In what follows, we shall freely interpret the results of 
\S\ref{6quadrics} as (Galois) bundle theoretic statements. The 
justification for this is that the Galois structure moves the fibres 
in just the right way for the global Clifford multiplications to make 
sense fibrewise. For example, if $U\subset V_F$ then consider the 
fibres $U_x$ and $(\tau U)_x$ at a point $x\in X$. By definition 
$$
(\tau U)_x = U_{\tau^{-1}(x)} \subset (V_F)_{\tau^{-1}(x)} = 
(S_F^+)_x,
$$
so that there is a well-defined multiplication map $U_x \otimes (\tau U)_x
\rightarrow (S_F^-)_{x}$.
\end{rem}

{\it Proof of \ref{mu=0}.}
Apply Lemma \ref{zero-norm} to $\sigma U \subset V_F$ with $g = \tau 
\sigma$. This says that the composition $\sigma U \otimes \tau U 
\hookrightarrow V_F \otimes S_F^+ \rightarrow S_F^-$ has rank 0 if 
$U$ is a line bundle and rank $\leq 1$ if $U$ has rank 2. Now apply 
$\tau$ and use the $S_3$ equivariance of the Chevalley algebra 
structure.
\qqed

We consider first a destabilising subbundle $E\subset V_F$ of rank 2.
Observe that if $\mu$ has rank 1 then  
the image of $\mu$ is a destabilising line subbundle; we postpone 
this case until after Proposition \ref{rank2case} and assume that $\mu =0$.
(The fact that the image of $\mu$ is isotropic
can be seen from the description of the Clifford multiplication map 
given in section \ref{2by2}.) 

The second thing to observe is that, with respect to the quadratic 
form,
$$
E^{\perp}/E \cong A_E \otimes B_E.
$$
This is because $E^{\perp}/E$ is an orthogonal bundle of rank 4 whose 
isotropic 2-planes in the generic fibre are parametrised by $A_E, 
B_E$ respectively.
(See \S\ref{6quadrics}).

\begin{prop}
\label{rank2case}
If $V_F$ has no destabilising line subbundle and
$E\subset V_F$ is a destabilising rank~2 subbundle then $F$ is 
S-equivalent to $F_{E}$ (\ref{split-form}), where $E$ is the 
pull-back of a stable rank~2 vector bundle on $X/\sigma$.
\end{prop} 

\prf
We already have the S-equivalence
$V_F \sim E\dual  \oplus E \oplus A_E \otimes B_E$; we examine the 
consequences of the vanishing of $\mu$.
Namely, $\tau^2 E = B_{\tau\sigma E}$ by Lemma \ref{AB2}, and so 
$E=B_{\tau^2\sigma E}$ by Lemma \ref{AB3}. So by Lemma \ref{AB1} we have
$$
A_E = A_{B_{\tau^2\sigma E}} = \tau^2\sigma E.
$$
Similarly $\sigma E = \sigma B_{\tau^2\sigma E} = A_{\sigma\tau^2\sigma 
E}= A_{\tau E}$ and so
$$
B_{\sigma E} = B_{A_{\tau E}} = \tau E.
$$

Next, $E\dual \cong 
\Delta^{-1}\otimes E$ where $\Delta = \det E$, and so  
we have
$$
\begin{array}{rcl}
V_F &\sim & (\oo \oplus \Delta^{-1})E \oplus A_E B_E,\\
S^+_F &\sim & (\oo \oplus \Delta^{-1})A_E \oplus B_E E,\\
S^-_F &\sim & (\oo \oplus \Delta^{-1})B_E \oplus E\ A_E.\\
\end{array}
$$
We now impose
the Galois condition. $S^+_F = \tau V_F$ implies that $\Delta $ is 
$\tau$-invariant ($\oo \oplus \Delta^{-1}$ is the only split rank 2
factor since $E,A_E,B_E$ are all stable). Similarly $\Delta$ is 
$\sigma$-invariant since $V_F = \sigma V_F$; it therefore descends to a line
bundle on $\bp^1$, and since it has degree 0 this forces $\Delta = \oo_X$.

Similarly $\sigma E \cong E$, and it follows from the above 
$A_E = \tau^2 E$, $B_E = 
\tau E$. As before, we note that since the $-1$-eigenspaces at fixed points
of $\sigma$ lie in the fibres of $A_E B_E$, the involution must act trivially
in the fibres of $E$, which therefore descends.
\qqed

The case of a destabilising line subbundle, to which we turn now, 
involves a bit more work.

\begin{prop}
\label{prop.1}
Suppose that a line subbundle $L\subset V_F$ is 
isotropic of degree 0, and that the Clifford multiplication $\tau L\otimes 
\tau^2 L \rightarrow V_F$ vanishes. Then $F \sim F_E$ for 
$E = L\oplus L^{-1}$.
\end{prop}

\prf
By (\ref{prod=0}) in \S\ref{6quadrics} and the hypothesis of the 
proposition, we have an isotropic subbundle $U = A_{\tau L} \cap 
B_{\tau^2 L}\subset V_F$ of rank 3. Using Lemma \ref{AB3}, moreover, 
one sees that $L\subset U$. Accordingly we have an S-equivalence
$$
V_{F} \sim L\oplus L\dual \oplus U/L \oplus (U/L)\dual \oplus 
U^{\perp}/U.
$$
We shall identify the last term first; note that $U^{\perp} = A_{\tau L} + 
B_{\tau^2 L}$, so we have to compute the rank 1 quotients $A_{\tau 
L}/U$ and $B_{\tau^2 L}/U$.

From (\ref{ABsequences}) it follows 
that there is a short exact sequence of vector bundles
$$
\ses{\tau L \otimes A_{\tau L}}{\tau L \otimes V_F}{B_{\tau L}}.
$$
Twisting by $\tau L^{-1}$ and using the fact that $\tau L \otimes 
\tau^2 L \rightarrow V_{F}$ is zero we get short exact sequences (in 
which the vertical arrows are inclusions):
$$
\begin{array}{rcccccl}
	0\rightarrow & A_{\tau L} & \rightarrow & V_F & \rightarrow & \tau 
	L^{-1}\otimes B_{\tau L} & \rightarrow 0  \\
	 & \uparrow &  & \uparrow &  & \uparrow &   \\
	0\rightarrow & U & \rightarrow & B_{\tau^2 L} & \rightarrow & \tau 
	L^{-1}\otimes \tau^2 L & \rightarrow 0,  \\
\end{array}
$$
and similarly
$$
\begin{array}{rcccccl}
	0\rightarrow & B_{\tau^2 L} & \rightarrow & V_F & \rightarrow & 
	\tau^2 
	L^{-1}\otimes A_{\tau^2 L} & \rightarrow 0  \\
	 & \uparrow &  & \uparrow &  & \uparrow &   \\
	0\rightarrow & U & \rightarrow & A_{\tau L} & \rightarrow & \tau^2 
	L^{-1}\otimes \tau L & \rightarrow 0.  \\
\end{array}
$$
We deduce that 
$$
V_F \sim (L  \oplus \tau^2 L^{-1}\otimes \tau L \oplus U/L) 
\oplus (\hbox{same bundle})\dual.
$$
To impose the condition that $F$ is a Galois bundle we shall compute 
directly the spinor bundles from this expression and (\ref{spinors}) 
in \S\ref{group}. First note that $\det A_{\tau L} = \det B_{\tau 
L} = \tau L^2$, so from the exact sequences above $\det U = \tau 
L \otimes \tau ^2 L$ and hence
$$
\det (U/L) = L^{-1}\tau L \tau ^2 L.
$$
Writing $N = L  \oplus \tau L\tau^2 L^{-1} \oplus U/L$, the spinor 
bundles are, up to a line bundle, $S_F^+ = \bigwedge^{\rm even} N$ and 
$S_F^- = \bigwedge^{\rm 
odd} N$. We find that to get trivial determinant one must twist by 
$\tau L^{-1}$, and the result is then (writing $E = L\otimes (U/L)$ 
and using the fact that (since it has rank 2) $(U/L)\dual = \tau L^{-1} \tau^2 L^{-1} \otimes 
E$):
$$
\begin{array}{rcl}
	V_F & \sim & L\oplus L^{-1} \oplus \tau L \tau^2 L^{-1}\oplus  
	\tau L^{-1} \tau^2 L \oplus U/L \oplus (U/L)\dual\\
	 & = &  L\oplus L^{-1} \oplus \tau L \tau^2 L^{-1}\oplus  
	\tau L^{-1} \tau^2 L \oplus \left( L^{-1}\otimes E\right) \oplus  
	\left( \tau L^{-1} \tau^2 L^{-1} \otimes E\right),\\
&&\\
	S_F^{+} & \sim &\textstyle  \bigwedge^{\rm even} N \\
	 & = & \tau  L\oplus \tau L^{-1} \oplus \tau^2 (L)  L^{-1}\oplus  
	\tau^2 (L^{-1})  L \oplus \left( \tau L^{-1}\otimes E\right) 
	\oplus \left( L^{-1}\tau^2 
	L^{-1} \otimes E \right), \\
&&\\
	S_F^- & \sim & \textstyle  \bigwedge^{\rm odd} N  \\
	 & = & \tau^2  L\oplus \tau^2 L^{-1} \oplus L \tau L^{-1}\oplus  
	 L^{-1} \tau L \oplus \left( \tau^2 L^{-1}\otimes E\right)  
	 \oplus \left( L^{-1}\tau L^{-1}
 \otimes E \right). \\ 
\end{array}
$$
Since $F$ is a Galois bundle these three vector bundles are permuted 
cyclically by $\tau$, and it follows that the bundle $E = L \otimes 
U/L$ is $\tau$-invariant. On the other hand, $E$ is also 
$\sigma$-invariant (since $U$ is), and $\sigma$ acts as the identity 
in the fibres of $E$ at fixed points since the $-1$-eigenspace in the 
fibre of $V_F$ is 1-dimensional and nonisotropic. It follows that $E$ 
descends to a bundle on $\bp^1$ (in fact the trivial bundle, by 
semistability of $V_F$, and can be written (notice that both summands 
are trivial)
$$
E = \oo \oplus L\tau L \tau^2 L.
$$
It follows directly that 
$$
V_F \sim (L \oplus L^{-1} \oplus \tau L \tau^2 L \oplus \tau L \tau^2 
L^{-1})\oplus (\hbox{same bundle})\dual
$$
and that $F$ has the desired form.
\qqed

We next show that if $V_F$ has a destabilising line subbundle
then this can be chosen with the property of Proposition \ref{prop.1}.

\begin{lemm}
Suppose that a destabilising line subbundle $L\subset V_F$ is polar,
with respect to the quadratic form on $V_F$, to its image 
$\sigma L \subset V_F$.
Then the Clifford multiplication $\tau L\otimes 
\tau^2 L \rightarrow V_F$ vanishes.
\end{lemm}

\prf
If $\sigma L = L$ then the result follows from Lemma \ref{mu=0}, so 
we may assume that $L, \sigma L$ are distinct subbundles. In this case 
the bundle map $L\oplus \sigma L \rightarrow V_F$ is injective by 
semistability (the image is isotropic since we assume the summands 
are polar). In other words $L$ and $\sigma L$ are distinct in all 
fibres.

If the multiplication $\tau L \otimes \tau^2 L \rightarrow V_F$ is 
nonzero then it must be injective (again by semistability, since the 
image is isotropic). This implies that $A_{\tau L} \cap B_{\tau^2 L} 
\subset V_F$ has rank 1 (using (\ref{prod=0})). On the other hand it 
contains $\sigma L$ (by Lemma \ref{mu=0}), and hence
$$
\sigma L = A_{\tau L} \cap B_{\tau^2 L}.
$$

We now consider the isotropic subbundles
$$
U = A_{\tau \sigma L} \cap B_{\tau^2 L},
\qquad
E = A_{\tau \sigma L} \cap A_{\tau L}.
$$
By Lemma \ref{mu=0} $U$ has rank 3. By \S\ref{6quadrics} 
(\ref{polarity}), polarity of $L$ and $\sigma L$ we see also that 
$E = \tau (A_{\sigma L} \cap A_{L})$ has rank 2. But $E,U$ are both 
contained in the rank 4 bundle $A_{\tau \sigma L}$, so
$$
N = U\cap E = A_{\tau \sigma L}\cap A_{\tau L} \cap B_{\tau^2 L}
$$
is an isotropic line subbundle. ($E$ is not contained in $U$ since 
then $A_{\tau L} \cap B_{\tau^2 L}$ would be $\geq 2$.) But $N\subset 
A_{\tau L} \cap B_{\tau^2 L} = \sigma L$, so in particular $\sigma L 
\subset  A_{\tau \sigma L}$. Applying $\sigma$ (and using Lemma 
\ref{AB3}) this shows that $L\subset B_{\tau^2 L}$. Applying $\tau^2$ 
then shows that $\tau^2 L\subset B_{\tau L}$, i.e. that $\tau L 
\otimes \tau^2 L \rightarrow V_F$ is the zero map, a contradiction.
\qqed

\begin{prop}
\label{prop.2}
If $V_F$ admits a destabilising line subbundle then it has such a line 
subbundle $L \subset V_F$ for which $\tau L\otimes 
\tau^2 L \rightarrow V_F$ vanishes.
\end{prop}

\prf
We assume that $L\subset V_F$ is a destabilising line subbundle 
without this property, i.e. for which $\tau L\otimes 
\tau^2 L \rightarrow V_F$ is injective; and by the previous lemma we 
can assume that $L$ is nowhere polar to $\sigma L$, i.e. that the 
homomorphism $L \otimes \sigma L \rightarrow \oo$ induced by the 
quadratic form on $V_F$ is an isomorphism. By \S\ref{6quadrics} 
(\ref{polarity}) it follows that in this case $V_F = A_{\tau \sigma 
L} \oplus A_{\tau L}$. We let $U=A_{\tau \sigma 
L} \cap B_{\tau^2 L}$, and we have a configuration of subbundles:
 \begin{center}
 \leavevmode
\epsfxsize=2in\epsfbox{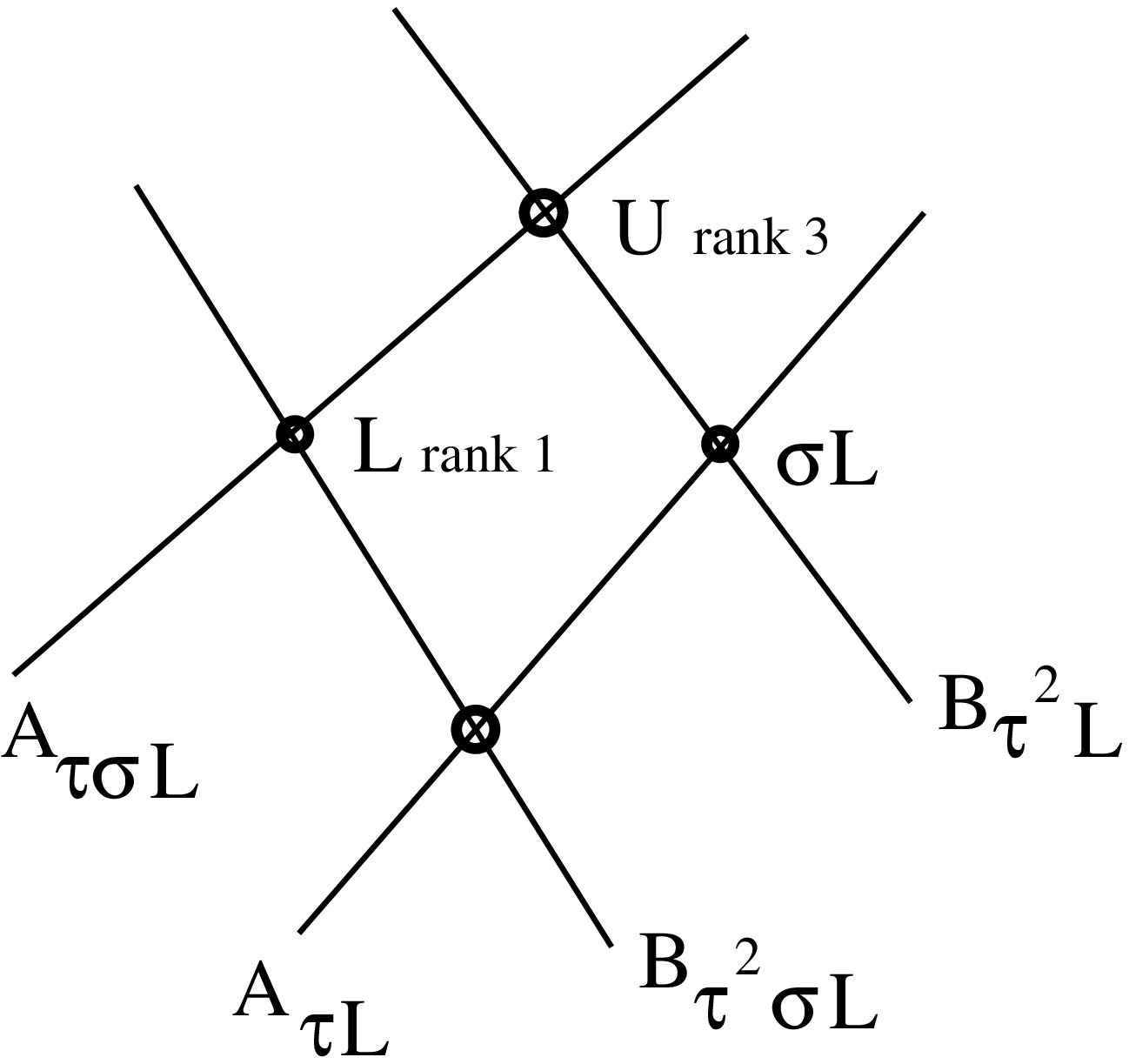}
\end{center}
Note in particular that 
$$
A_{\tau L} \cap B_{\tau^2 L} = \tau L\tau^2 L = \sigma L \cong L^{-1}.
$$
As in the proof of Proposition \ref{prop.1} we have exact sequences 
(the vertical arrows are injective)
$$
\begin{array}{rcccccl}
	0\rightarrow & A_{\tau \sigma L} & \rightarrow & V_F & \rightarrow & \tau 
	\sigma L^{-1}\otimes B_{\tau \sigma L} & \rightarrow 0  \\
	 & \uparrow &  & \uparrow &  & \uparrow &   \\
	0\rightarrow & U & \rightarrow & B_{\tau^2 L} & \rightarrow & \tau 
	\sigma L^{-1}\otimes \tau^2 L & \rightarrow 0.  \\
\end{array}
$$
From this we see that $\det U = \tau^2 L \tau L^{-1}$ and
$$
\textstyle
\bigwedge^2 U = U\dual \otimes \det U = U\dual \otimes \tau^2 L \tau L^{-1}.
$$
Using these facts and \S\ref{group} (\ref{spinors}) we deduce that 
$$
\begin{array}{rcl}
	V_F & = & (L\oplus U) \oplus (L\oplus U)\dual  \\
	S_F^+ & = & (\tau L \oplus L\tau L \otimes U) \oplus (\tau L \oplus L\tau 
	L \otimes U)\dual\\
	S_F^- & = & (\tau^2 L \oplus \tau L \otimes U ) \oplus 
	(\tau^2 L \oplus \tau L \otimes U )\dual.
\end{array}
$$
We now impose the condition that $F$ is a Galois bundle to deduce 
that 
$$
\tau U = L \tau L \otimes U. 
$$
Also $\sigma U = U$ since $\sigma$ 
interchanges $A_{\tau \sigma L}$ and $B_{\tau^2 L}$ (by Lemma 
\ref{AB3}). At fixed points $\sigma$ acts trivially in the fibres of 
$U$ since the latter is isotropic, and hence the projective bundle 
$\bp(U)$ is $S_3$-invariant and descends to a $\bp^2$-bundle on 
$\bp^1$. This is necessarily the projectivisation of a split vector 
bundle, and hence 
$U = P\oplus Q \oplus R$ for some line bundles $P,Q,R$ on $X$. Note 
that these line bundles need not themselves descend to $\bp^1$ may be 
permuted and twisted by the action of $S_3$. However, since $U$ is 
$\sigma$-invariant at least one summand must also be 
$\sigma$-invariant, and this summand will be a destabilising subbundle 
of $V_F$ with the desired properties.
\qqed

{\it Proof of Theorem \ref{ssboundary}.}
By Proposition \ref{rank2case} we are reduced to the case where 
$V_F$ has a destabilising line subbundle 
$L\subset V_F$. By Proposition \ref{prop.2} we can assume that Clifford
multiplication on $\tau L \otimes \tau^2 L$ is zero, and by \ref{prop.1}
this implies that $V_F \sim \c^2 E \oplus \tau E \tau^2 E$ with $E = L \oplus 
L^{-1}$. We just have to check that $E$ descends to a bundle on $X/\sigma$.
On the other hand the Galois condition tells us that $V_F \sim \c^2 \sigma E 
\oplus \tau^2 \sigma E \tau \sigma E$, and this implies that $\sigma E \cong 
E$. At a fixed point $x\in X$ of $\sigma$ the 1-dimensional $-1$-eigenspace 
of $\sigma$ in the fibre of $V_F$ is $\bigwedge^2 (\tau E_x)$ (since $\tau E$ 
and $\tau^2 E$ are interchanged by $\sigma$), and so in particular $\sigma$ 
acts trivially in the fibre $E_x$. This shows that $E$ is pulled back from a 
bundle on $X/\sigma$, as required.
\qqed

\section{Triality and trigonality}
\label{galois3}

In this section (in \S\ref{trig-gal-bun}) 
we shall introduce the main object of the paper, 
which is the moduli space $\nn_C$ of Galois $\spin(8)$-bundles on the Galois 
closure of a trigonal curve $C\map{3:1} \bp^1$. Unless the curve is 
cyclic this Galois curve is a connected double cover of $C$, and we are 
interested in the Galois bundles on this double cover. The moduli 
space $\nn_C$ is an extension of the variety $\su(2)$, which contains its 
singular locus. 

\subsection{Galois trigonal curves}
\label{trig-pics}

Let $C$ be a nonhyperelliptic curve of genus $g\geq 3$ with a fixed 
trigonal pencil $g^{1}_{3}$. We shall denote by $\bb 
\subset \bp^{1}$ the branch divisor; by Riemann-Hurwitz this has 
degree 
$2g+4$. We shall suppose that $\bb$ contains $\delta$ points of
multiplicity two, i.e. at which the cover has cyclic triple 
branching. 
There are then $2g+4 -2\delta$ simple branch points:
\begin{center}
 \leavevmode
\epsfxsize=2.5in\epsfbox{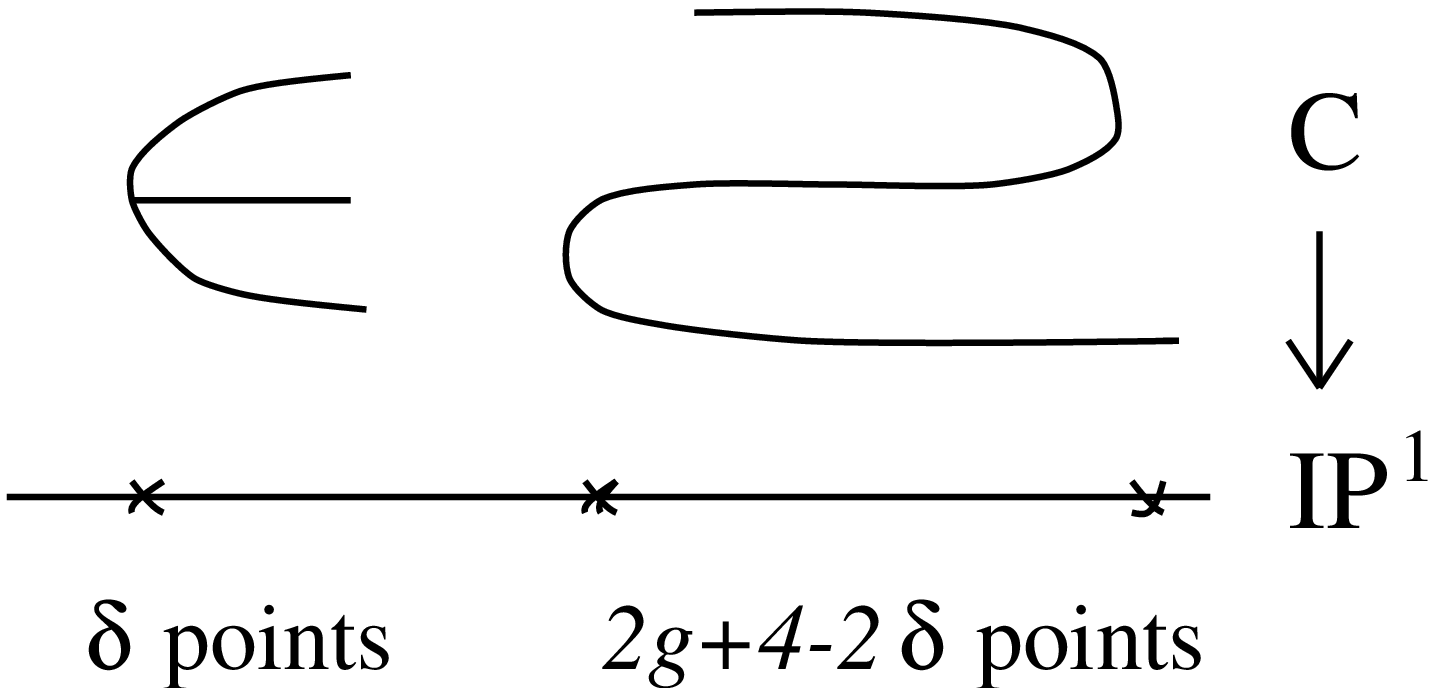}
\end{center}
We shall need to distinguish the case of a cyclic cover $C\map{3:1} 
\bp^{1}$, in which case $\delta = g+2$. (We shall see in a moment 
that 
this is also a sufficient condition.)

We construct a Galois cover $\gal(C) \rightarrow \bp^{1}$ as 
follows. $\gal(C)\subset C\times C$ is the closure of the set 
$\{(q,r)\ |\ q\not= r,\ p+q+r \in g^{1}_{3}\ \hbox{for some }p\in 
C\}$. This is a 
double cover $\gal(C) \rightarrow C$ by $(q,r) \mapsto p= 
g^{1}_{3}-q-r$, 
and is smooth with branching behaviour:
\begin{center}
 \leavevmode
\epsfxsize=2.5in\epsfbox{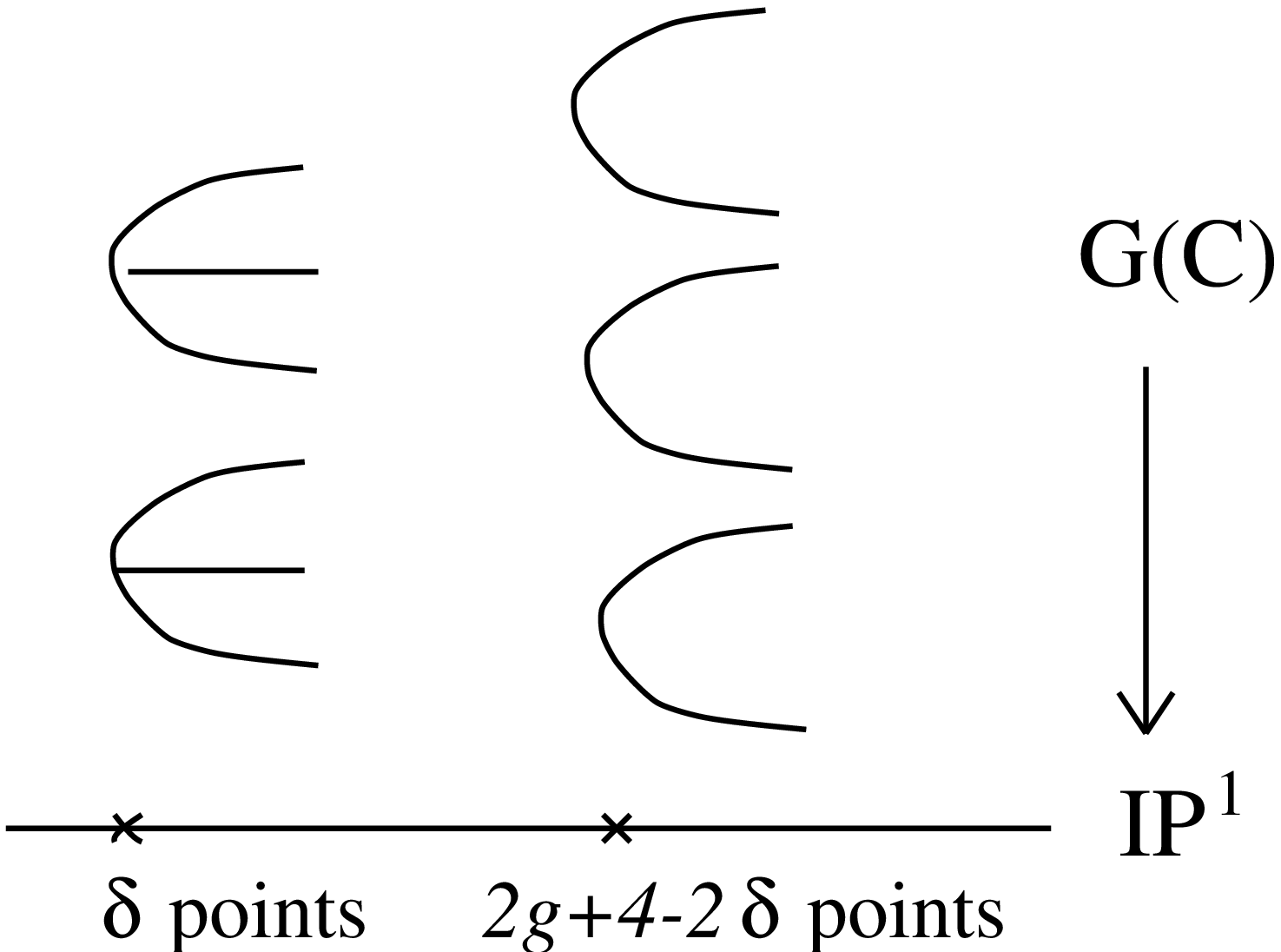}
\end{center}

Notice that over each divisor $p+q+r \in g^{1}_{3}$ the points of 
$\gal(C)$ 
correspond to the orderings of $p,q,r \in C$. Thus the permutation 
group $S_{3}$ acts on $\gal(C)$ on the right; in particular $C$ is 
the 
quotient of $\gal(C)$ by the 2-cycle $\sigma = (23)$. 
In general the general case 
$S_3$ is the Galois group of $\gal (C)$; however, in the cyclic 
case $\gal (C)$ breaks up into two connected components, each a copy
of $C$ acted on by $\z/3$.

Restricting for a moment to the Galois $S_3$ case, we find 
by Riemann-Hurwitz that $\gal(C)$ has genus $3g+1-\delta$. 
Over simple branch points in $\bb \subset \bp^{1}$ the
three 
ramification points  
are respective fixed points of the three 2-cycles in $S_{3}$.
Over each double branch point there are two fixed points of the cyclic
subgroup $\<\tau\>\cong \z/3$, at which $\tau$ acts in the tangent 
space
as $\omega$, $\omega^2$ respectively, where $\omega = e^{2\pi i/3}$ 
is a 
primitive cube root of unity.

Consider now
the quotient $H' = \gal(C)/\<\tau\>$. This is a hyperelliptic curve
branched over 
$\bb \subset \bp^{1}$ and singular over the $\delta$ double points of 
$\bb$. 
By Riemann-Hurwitz (again, we are assuming $\gal (C)$ is a connected 
Galois 
$S_3$ curve)
its normalisation $H$ has Euler characteristic $e(H) = 2\delta -2g$. 
If
$\delta = g+2$ then $e(H) =4$ so $H$ is a union of two lines and 
$\gal (C)$ 
is a pair of copies of $C$, a contradiction. We have therefore shown:

\begin{lemm}
$C$ is cyclic over $\bp^1$ if and only if $\delta = g+2$.
\end{lemm}

If $\delta \leq g+1$, on the other hand, then $g(H) = g+1 - \delta$.
We can summarise the situation in the following diagram. 

\begin{equation}
\label{Gal(C)-diagram}
	\begin{array}{ccccc}
	 &  & ^{3g+1-\delta}\gal(C)=X &  &   \\
	
	 & \swarrow  &  & \searrow\pi &   \\
	^{g+1-\delta}H' &  &  &  &  ^{g}C =Y \\
	 
	 & \searrow &  & \swarrow g^{1}_{3} &   \\
	
	 &  & \bp^{1} &  &   \\
\end{array}
	\label{S3galois}
\end{equation}

In the rest of the paper we shall apply the discussion of the 
previous 
sections to $X = \gal(C)$ and $Y=C$, in the general trigonal case,
and to $X=C$ when this is a cyclic trigonal curve. Whenever we 
refer to $\gal(C)$ we shall always assume that we are in the general 
case.

\subsection{Trigonal Galois bundles}
\label{trig-gal-bun}

We now turn to the moduli space we are interested in. We assume that 
$C$ is a non-cyclic trigonal curve and we let 
$\nn_{C} := \nn_{\gal(C)} \subset \mm_{\gal(C)}(\spin(8))$.
(Though see Remark \ref{dimcount}(i) regarding the cyclic case.)

\begin{rem}\rm
If $F$ is such a Galois bundle then the sheet-involution of
$\gal(C)$ over $C$, corresponding to $\sigma \in S_3$, lifts to the 
orthogonal bundle $V_F$; applying the 3-cycle $\tau \in S_3$ one 
sees that
the fibres of $V_F$ over any $p+q+r \in g^1_3$ can be 
identified 
(though not canonically) with the three representations $V,S^+,S^-$.
If $F$ is semistable but not stable then $F = F_E$ for some 
$E \in \ss\uu_C(2)$ (by Theorem \ref{ssboundary}).
The following table shows explicitly how, in this case, the 
fibres of 
$V_{F}, S_{F}^{\pm}$ are constructed from those of the rank 2 bundle 
$E$ over 
a divisor $p+q+r \in g^{1}_{3}$. We abbreviate 
$\c^2 \otimes E_{p} \oplus E_{q}\otimes E_{r}$ to $2E_p+E_q E_r$:
$$
\begin{array}{c|c||c|c|c}
	
	C & \gal(C)  & V_{F} & S^{+}_{F} & S^{-}_{F}  \\
	
	\hline
	
	p & (q,r) &      2E_p+E_q E_r     & 2E_q+E_r E_p &  2E_r+E_p E_q \\
	 & (r,q) &     2E_p+E_r E_q  & 2E_r+E_q E_p &  2E_q+E_p E_r \\
	
\end{array}
$$
\end{rem}
\noindent
The following is the main result of this paper:

\begin{theo}
\label{dimtheorem}
For any non-cyclic trigonal curve $C$ the moduli space
$\nn_C$ is smooth of dimension $7g-14$ away from its semistable 
boundary $\ss\uu_C(2)$.
\end{theo}

\begin{rems}\rm
\label{dimcount}
\itemitem{(i)}
We have already remarked in \S\ref{hypcase} that most of our constructions
work also for cyclic trigonal curves. to be precise,
if $C$ is cyclic we 
define $\nn_{C}= \ff_C \subset \mm_{C}(\spin(8))$, the fixed-point-set
of the conjugation action of $\z/3$ after choosing an embedding $\z/3
\hookrightarrow S_3$. In this case both the smoothness (at stable points)
and the dimension statements
in Theorem \ref{dimtheorem} remain valid (see Remark \ref{cyclic-dim} below);
what is missing is the identification of the semistable boundary with 
$\ss\uu_C(2)$ (though it is tempting to expect that this also remains true).
\itemitem{(ii)}
Notice that the moduli space of Galois bundles on the 
hyperelliptic curve $H$ of diagram (\ref{Gal(C)-diagram}) 
has dimension $7g-14 - 7\delta$ 
(see \S\ref{hypcase}); when $\delta=0$ and
$H'= H$ is smooth this number is $7g-14$.
This raises the question of the relationship between these two moduli spaces.
Are they isomorphic, for example, or birational?
\end{rems}

In view of Proposition \ref{tgt-space}, Theorem \ref{ssboundary} and 
the fact that in the trigonal situation the quotient $X/\sigma$ is 
just the curve $C$ itself,
all that remains to prove of Theorem \ref{dimtheorem} is 
the dimension statement; this will be 
carried out next.

\subsection{Dimension calculation}
\label{gen-dim}

By Proposition \ref{tgt-space} we have to compute the dimension of the invariant subspace 
$H^1(\gal(C),\ad F)^{S_3}$ (or in the cyclic case $H^1(C,\ad 
F)^{\z/3}$). We shall do a little more, in fact, and compute the 
decomposition of $H^1(\gal(C),\ad F)$ into irreducible 
representations under this group action. In the same way, of course, 
we could as 
easily compute $T_{F}\nn_X$ at a stable Galois bundle for any $X$, but 
we shall stick to the trigonal situation.

We shall consider the $S_3$ case first and return to $\z/3$ in a 
moment. 
For use here and below  
we recall the following result of \cite{AB} Theorem 4.12.

\begin{atiyah-bott}
\label{at-bott}
Suppose that $\gamma: X \rightarrow X$ is an automorphism of a 
compact 
complex manifold $X$ with a finite set $\fix(\gamma)\subset X$ of 
fixed 
points, 
and suppose that $\gamma$ lifts to a holomorphic vector bundle 
$E\rightarrow X$. 
Then 
$$
\sum (-1)^p \tr \gamma|_{H^p(X,E)} = \sum_{x\in \fix(\gamma)} {\tr 
\gamma|_{E_x}
\over \det(1-d\gamma_x)}.
$$ 
\end{atiyah-bott}

Suppose now that 
$F\in \nn_C$ is a stable Galois $\spin(8)$-bundle for a general 
trigonal curve $C$.
We use Atiyah-Bott to compute 
the trace of the group elements $\sigma,\tau \in S_3$ acting on 
$H^1(\gal(C), \ad F)$. By stability $H^0(\gal(C), \ad F)=0$ 
and so for nontrivial 
$\gamma \in S_3$:
$$
\tr \gamma|_{H^1(\gal(C), \ad F)} = - \sum_{x\in \fix(\gamma)} {\tr 
\gamma|_{\ad F_x}
\over \det(1-d\gamma_x)}.
$$

\underline{$\gamma = \sigma$}:
$\sigma$ has $2g+4-2\delta$ fixed points, and at 
each of these $d\sigma_x = -1$. On $\ad F_x$ 
we need 
to compute $\tr \sigma$. But the bundle $F$ is Galois, so 
this is just the trace of $\sigma$ in the triality action on $\so_{8}$,
which by section \ref{liealg} is 14. Thus
$$
\tr \sigma |_{H^1(\gal(C), \ad F)} = -7(2g+4-2\delta).
$$

\underline{$\gamma = \tau$}:
we have seen that $\tau$ has $\delta$ {\it pairs} of fixed points, at 
which in the tangent space it acts as $\omega$ and $\omega^2$ 
respectively. From section \ref{liealg}, on the other hand, $\tr \tau 
=7$ in the adjoint representation. So by Atiyah-Bott 
(and using $1+\omega +\omega^2 =0$) we get
$$
\tr \tau |_{H^1(\gal(C), \ad F)} = -7\delta \left\{
{1\over 1-\omega }+{1\over 1-\omega^2}
\right\} = -7\delta.
$$

Note that $\dim H^1(\gal(C), \ad F) = 28(3g-\delta)$ (since $\gal(C)$ 
has genus $3g+1-\delta$); and so we have computed the character of 
$H^1(\gal(C), \ad F)$ under
the $S_3$ action. 

Recall that $S_3$ has three irreducible representations $\c$, 
$\c_{\varepsilon}$, $\c^2_{\rho}$ and character table:

\begin{equation}
\label{S3-table}
\begin{array}{l|rrr}
&1&\sigma&\tau\\
\hline
\one&1&1&1\\
\varepsilon&1&-1&1\\
\rho&2&0&-1\\
\end{array}
\end{equation}
Here $\c = \c_{\one}$ is the trivial representation, 
$\c_{\varepsilon}$ the 
sign representation
and $\c^2_{\rho}$ the 2-dimensional reflection representation.

Denote the character of $H^1(\gal(C), \ad F)$ by $a\one +b\varepsilon + 
c\rho$. Comparing with the character table gives:
$$
\begin{array}{lcl}
a+b+2c &=& 28(3g-\delta)\\
a-b &=& -7(2g+4-2\delta)\\
a+b-c &=& -7\delta .\\
\end{array}
$$
Solving these equations yields 
$a= 7g-14$, $b=21g+14 -14\delta$, $c= 28g-7\delta$,
and hence an $S_3$-decomposition:
$$
H^1(\gal(C), \ad F) = \c^{7g-14} \oplus (21g+14-14\delta) \c_{\varepsilon} 
\oplus 
(28g-7\delta) \c^2_{\rho}.
$$

\begin{rem}\rm
\label{cyclic-dim}
One can make an analogous calculation in the case of a cyclic trigonal curve.
Here we apply the same method as above, with group $\<  \tau \> \cong \z/3$ 
acting
on $C$ itself. 
Pick a stable Galois $\spin(8)$-bundle $F\in \nn_C$.
As above we compute 
the trace of the group elements $\tau ,\tau^2$ acting on 
$H^1(C, \ad F)$, where by stability $H^0(C, \ad F)=0$. 

By Atiyah-Bott:
$$
\begin{array}{rcl}
\tr \tau |_{H^1(C, \ad F)} 	 & = & -7(g+2)/(1-\omega),  \\
\tr \tau^2 |_{H^1(C, \ad F)}	 & = & -7(g+2)/(1-\omega^2).\\
\end{array}
$$
Suppose that the eigenspace decomposition under $\tau$ is
$H^1(C, \ad F) = \c^a \oplus \c^b \oplus \c^c$ for eigenvalues 
$1,\omega,\omega^2$ respectively. Then we have to solve the equations
$$
\begin{array}{lcl}
a+b+c &=& 28(g-1)\\
a+\omega b +\omega^2 c&=& -7(g+2)/(1-\omega)\\
a+\omega^2 b +\omega c &=& -7(g+2)/(1-\omega^2) .\\
\end{array}
$$
Adding the three equations gives again $a= 7g-14$.
\end{rem}

\section{A nonabelian Schottky configuration}
\label{recillas}

We shall show in this section that the moduli space $\nn_C$ contains a 
nonabelian `fattening' of the classical Schottky configuration for 
trigonal curves. The Schottky configuration is the union of Prym 
Kummer varieties embedded in $\su(2)$; when $C$ is trigonal two things 
happen: (a) each Prym is, by the Recillas correspondence, 
the Jacobian of a tetragonal curve 
$R_{\eta}$, and (b) $\su(2)$ is by Theorem \ref{dimtheorem} 
the `singular locus' of a bigger moduli 
space $\nn_C$. We show that $\nn_C$ contains a configuration of the 
moduli varieties $\ss\uu_{R_{\eta}}(2)$---each singular, of course, 
along the corresponding Prym Kummer.

\subsection{The Schottky configuration}

For each nonzero 2-torsion point $\eta \in J_C[2] \backslash
\{\oo\}$ we have an associated unramified double cover
$$f : C_{\eta} \rightarrow C.$$ 
We shall denote by $\varphi$ the involution of $C_{\eta}$ given by 
sheet-interchange over~$C$; it will denote 
also the induced involution of $\pic(C_{\eta})$.
The kernel of the norm map on divisors has two isomorphic connected 
components:
$$
\nm^{-1}(\oo_C) = P_{\eta} \cup P_{\eta}^-,
$$
where $P_{\eta} = (1 - \varphi)J^0_{C_{\eta}}$ and $P_{\eta}^- = (1 - 
\varphi)J^1_{C_{\eta}}$. 

We shall also consider the 2-component subvariety
$\nm^{-1}(\eta) \subset J_{C_{\eta}}$; this is a torsor over $\nm^{-1}(\oo_C)$
whose components are exchanged under translation by elements of $P_{\eta}^-$. 
Note that since $\ker f^* = \{\oo,\eta\}$ and $f^* \circ \nm = 1+\varphi$
it follows that the set of anti-invariant (degree 0) line bundles on $C_{\eta}$ 
is
$$
\ker (1+\varphi) = \nm^{-1}(\oo_C)\cup \nm^{-1}(\eta).
$$
$P_{\eta} \subset \nm^{-1}(\oo_C)$ is distinguished as the component 
containing the origin $\oo\in J_{C_{\eta}}$; the component 
$P_{\eta}^- \subset \nm^{-1}(\oo_C)$ contains the Prym-canonical curve
$C_{\eta} \map{1 - \varphi}P_{\eta}^-$. The components of $ \nm^{-1}(\eta)$, on the other hand, are indistinguishable.

\begin{rem}\rm
The $[-1]$-action on $ J_{C_{\eta}}$ preserves each of these four 
components (acting as sheet-involution $\varphi$): clearly it permutes
the components isomorphically, and so it suffices to note that each 
component contains fixed points, i.e. points of $J_{C_{\eta}}[2]$. 
\end{rem}

For each $\eta \in J_C[2] \backslash \{\oo\}$ there exist canonical maps
$$
\nm^{-1}(\eta) \rightarrow \su(2) \rightarrow |2\Theta|
$$
where the first map is just direct image $x \mapsto f_* x$; this gives a 
semistable rank 2 vector bundle of determinant $\det f_* x = \nm(x)\otimes 
\det f_* \oo_{C_{\eta}} = \oo_C$. This map obviously factors through 
the sheet-involution, and so the image is a pair of Prym-Kummer varieties.
As $\eta \in J_C[2]$ varies,
the configuration of these Kummers in $\su(2) \subset |2\Theta |$, including 
(for $\eta = \oo$) the Jacobian Kummer, is called the {\it Schottky 
configuration}. 

It is usually convenient to view the Kummer varieties of the Schottky 
configuration
as the images of $\nm^{-1}(\oo_C) = P_{\eta} \cup P_{\eta}^-$, though this
can only be done up to the translation action of $J_C[2]/\eta$. For this
we introduce the $J_C[2]$-torsor
$$
S(\eta) = \{\zeta\ |\ \zeta^2 = \eta\} \subset J_C.
$$
Notice that pull-back gives an 
isomorphism of groups
$$
f^* : J_C[2]/\eta \idfy \nm^{-1}(\oo_C)[2] := \nm^{-1}(\oo_C)\cap 
J_{C_{\eta}}[2],
$$
and correspondingly an isomorphism of torsors
$$
f^* : S(\eta)/\eta \idfy \nm^{-1}(\eta)[2].
$$

Each choice of $\zeta \in S(\eta)/\eta$ determines an identification
$\nm^{-1}(\oo_C) \idfy \nm^{-1}(\eta)$ by translation $x\mapsto f^*(\zeta)
\otimes x$. This sends 2-torsion points to 2-torsion points, and determines
a map
$$
\sch_{\zeta}:
\peta \cup \peta^- \rightarrow \su(2) 
$$
by $x \mapsto f_*(f^*(\zeta)
\otimes x)=\zeta \otimes f_* x$. The image is 
independent of the choice of $\zeta$ and is precisely the fixed-point 
set of the involution $\otimes \eta$ of $\su(2)$. Moreover, the 
linear span of the image of $P_{\eta}$ in $|2\Theta|$
can be canonically identified with the linear 
series $|2\Xi|$ where $\Xi$ is the canonical theta divisor on the 
dual abelian variety, and represents the principal polarisation on 
$\peta$ induced from that on $J(C_{\eta})$. Thus for each $\zeta \in 
S(\eta)/\eta$
there is a commutative diagram
\begin{equation}
\label{prymkum}
\begin{array}{ccc}
\peta & \map{\sch_{\zeta}} & \su(2) \\
&&\\
\downarrow &&\downarrow\\
&&\\
|2\Xi| &\hookrightarrow & |2\Theta|.\\
\end{array}
\end{equation}

\subsection{$S_4$-curves over $\bp^1$}

When the curve $C$ is trigonal the Pryms of the Schottky configuration are
in fact Jacobians, by a well-known construction of Recillas.
We shall next review this construction, as formulated 
by Donagi \cite{D}. 

To any trigonal curve $C\map{3:1} \bp^{1}$ 
together 
with a nontrivial 2-torsion point $\eta \in J_{C}[2]$ one can associate in a 
natural way a tetragonal curve $R_{\eta}\map{4:1} \bp^{1}$ whose 
polarised Jacobian is isomorphic to the Prym variety 
of $C_{\eta} \rightarrow C$. 
By definition $R_{\eta}$ is one of two isomorphic connected 
components 
of the fibre product
$$
\begin{array}{ccc}
R_{\eta}\cup R_{\eta}' &\subset & S^3 C_{\eta}\\
&&\\
\downarrow &&\downarrow \\
&&\\
\bp^1 &\hookrightarrow & S^3 C \\
\end{array}
$$
where $\bp^1 \hookrightarrow  S^3 C$ is the trigonal pencil, and the 
right-hand vertical map is the 8:1 cover induced from 
$C_{\eta}\rightarrow C$.
Conversely $C_{\eta} = S^2_{\bp^1} R_{\eta} \subset S^2 R_{\eta}$ 
with 
the obvious involution. In particular, $C_{\eta}$ and $R_{\eta}$ are 
in 
(3,2)-correspondence, and we shall denote the incidence curve by 
$\ii \subset C_{\eta}\times R_{\eta}$.

We shall refer to $R_{\eta}$ as the {\it 
Recillas curve}, and we shall consider its Galois closure 
$\gal(R_{\eta})\rightarrow \bp^{1}$ with Galois group $S_{4}$. As in 
section \ref{trig-pics} the points of $\gal(R_{\eta})$ correspond to 
orderings of the $g^{1}_{4}$-divisors in $R_{\eta}$.
 
\begin{notat}\rm
\label{S4notation}
We shall view $S_{4}$ as the permutation group of $\{0,1,2,3\}$, 
containing $S_{3}$, the permutations of $\{1,2,3\}$, as a subgroup. 
We shall fix elements $\sigma = (23)$, $\tau = (123)$ (as in earlier 
sections) and also $\sigma' = (01)$, $\varphi = (0213)$. 
(Note that $\varphi$ is uniquely determined up to inverse by the 
condition
$\sigma \sigma'=\varphi^2$. The reason for using the same notation $\varphi$
as for the sheet-involution of $C_{\eta}$ will appear in a moment.)
We shall then refer to the subgroups `Klein', generated by the 
conjugates of 
$\sigma \sigma'=\varphi^2$, `dihedral' $= \<\sigma ,\varphi\>$ and 
$A_{4}$ the alternating subgroup.

$S_4$ has five irreducible representations 
and character table:
$$
\begin{array}{c|rrrrr} 
	& 1 & \sigma & \tau & \varphi & \sigma \sigma'  \\ 
	\hline
\one  & 1 & 1 & 1 & 1 & 1  \\	
    \varepsilon  & 1 & -1 & 1 & -1 & 1  \\
	\rho & 3 & 1 & 0 & -1 & -1  \\
	\rho \otimes \varepsilon & 3 & -1 & 0 & 1 & -1  \\
	\two & 2 & 0 & -1 & 0 & 2
\end{array}
$$
Here $\rho$ is the 3-dimensional reflection representation 
and $\two$ is the 2-dimensional reflection representation of 
$S_{4}/{\rm Klein} \cong S_{3}$.
\end{notat}

The various curves in the Recillas construction are related by the 
following 
diagram which extends (\ref{S3galois}):

\begin{equation}
\label{S4galois}
\begin{array}{ccccc}
	 &  & \gal(R_{\eta}) &  &   \\
	 & \swarrow &  & \searrow &   \\
	\gal(C)\times_{C}C_{\eta} &  &  &  & \ii  \\
	&&&&\\
	\downarrow & \searrow &  &\textstyle 2:1\swarrow &   \\
        &&&&\\
	\gal(C) &  & C_{\eta} &  & \downarrow\textstyle 3:1 \\
        &&&&\\
	\downarrow & \searrow & \downarrow &  &   \\
        &&&&\\
	^{g+1-\delta}H &  & ^{g}C &  & ^{g-1}R_{\eta}  \\
	 & \searrow & \downarrow & \swarrow &   \\
         &&&&\\
	 &  & \bp^{1} &  & \\
\end{array}
\end{equation}
By comparison of \S\ref{trig-pics} and \cite{D} p.74, one sees that 
the Galois cover 
$\gal(R_{\eta}) \rightarrow \bp^1$ has branching behaviour:
\begin{center}
 \leavevmode
\epsfxsize=2.5in\epsfbox{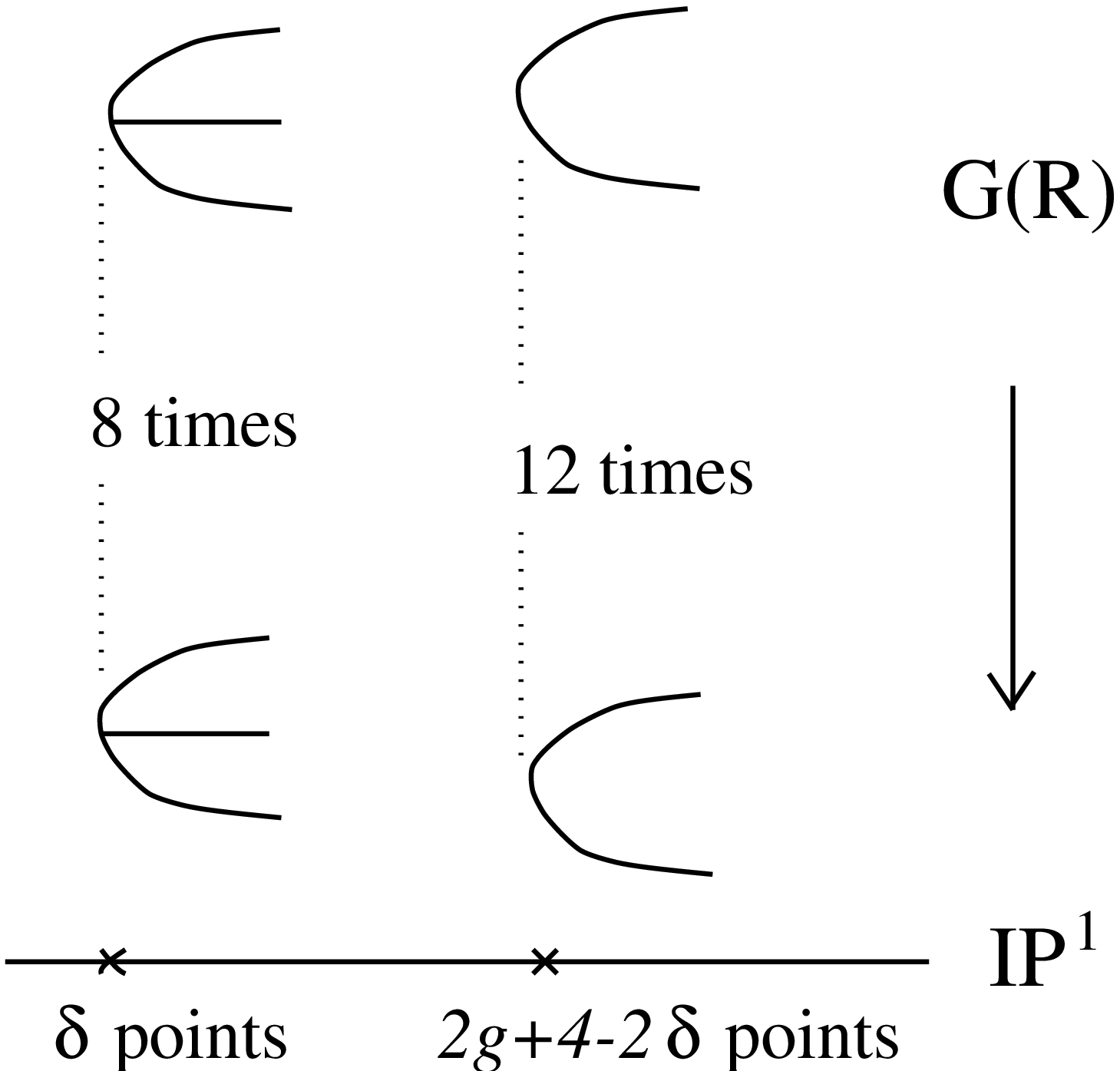}
\end{center}
By Riemann-Hurwitz it follows that $\gal(R_{\eta})$ 
has 
genus $12g+1-4\delta$, and in particular we see: 

\begin{lemm}
\label{etale}
The 4:1 map 
$\gal(R_{\eta})
\rightarrow \gal(C)$ is unramified. 
\end{lemm}

Diagram (\ref{S4galois}) is 
Galois-dual to the tower of subgroups (and these diagrams  
coincide with those of \cite{D}~p.73):
\begin{equation}
\label{S4galoisdual}
\begin{array}{ccccc}
	&  & 1 &  &   \\
	 & \swarrow &  & \searrow &   \\
         
	\<\sigma \sigma'\> &  &  &  &  \<\sigma \> \\
       
	\downarrow & \searrow &  & \swarrow &   \\
        &&&&\\
	{\rm Klein} &  & \<\sigma, \sigma'\> &  & \downarrow  \\
        &&&&\\
	\downarrow & \searrow & \downarrow &  &   \\
        &&&&\\
	A_{4} &  & {\rm dihedral} &  & S_{3}  \\
	 & \searrow & \downarrow & \swarrow &   \\
        &&&&\\
	 &  & S_{4} &  & 
\end{array}
\end{equation}
Note that the quotient dihedral$/\<\sigma, \sigma'\> \cong \z/2$ is generated
by $\varphi$, which therefore acts as sheet-interchange of $C_{\eta}$ over $C$,
consistently with our earlier notation.

\begin{prop}
\label{Hodgedecomp}
As $S_{4}$-modules we have
$$
H^{1}(\gal(R_{\eta}),\oo) = (g+1-\delta)\c_{\varepsilon} \oplus
g \c^{2}_{\two} \oplus 
(g-1) \c^{3}_{\rho} \oplus
(2g+1-\delta)\c^{3}_{\rho\otimes \varepsilon}.
$$
\end{prop}

\prf
One proceeds the same way as in the dimension computation of \S\ref{gen-dim}.
By Atiyah-Bott we have, for 2-cycles $\sigma \in S_4$,
$$
\tr \sigma|_{H^1(\gal(R_{\eta}),\oo)} =
1-{1\over 2}|\fix(\sigma)|=
-2g-3+2\delta;
$$
for 3-cycles $\tau \in S_4$,
$$
\tr \tau|_{H^1(\gal(R_{\eta}),\oo)} =
1-\delta \left\{ {1\over 1-\omega} +{1\over 1-\omega^2} \right\} 
=
1-\delta;
$$
while the trace of other nontrivial group elements is 1,
and the identity element of $S_4$ has trace $g(\gal(R_{\eta}))= 
12g+1-4\delta$.
If the representation on $H^1(\gal(R_{\eta}),\oo)$ has character
$a\one + b\varepsilon + c\rho + d\rho \otimes \varepsilon + e\two$
for integers $a,b,c,d,e$, then 
comparing with the character table \ref{S4notation} 
we have to solve:
$$
\begin{array}{rrrrrrr}
12g+1-4\delta &=&a&+b&+3c&+3d&+2e \\
-2g-3+2\delta&=&a&-b&+c&-d& \\
1-\delta&=&a&+b&&&-e \\
1&=&a&-b&-c&+d& \\
1&=&a&+b&-c&-d&+2e \\
\end{array}
$$
These equations give $(a,b,c,d,e)=(0,g+1-\delta,g-1,2g+1-\delta,g)$.
\qqed

We can now describe the polarised isomorphism $P_{\eta} \cong 
J_{R_{\eta}}$ in this setting. Both abelian varieties embed by 
pull-back in the Jacobian of the top curve $\gal(R_{\eta})$; and 
this Jacobian is acted on by the group-ring $\z[S_4]$. 

\begin{prop}
\label{recillasisom}
The action of $1+\varphi^2 \in \z[S_4]$ on the Jacobian 
$J_{\gal(R_{\eta})}$ restricts to the polarised isomorphism of 
Recillas
$\psi: J_{R_{\eta}}\idfy P_{\eta}$. Moreover, $2 \psi^{-1}:P_{\eta} 
\rightarrow J_{R_{\eta}}$ is the restriction of $1+\tau + \tau^2 \in 
\z[S_4]$.
\end{prop}

\prf
It is easy to identify the subspaces covering (the pull-backs of)  
these two abelian varieties in the $S_{4}$-decomposition of 
Proposition \ref{Hodgedecomp}: for $J_{R_{\eta}}$ we look for 
$S_3$-invariant
vectors, while for $P_{\eta}$ we look for $\<\sigma, \sigma'\>$ 
invariant 
vectors on which the sheet-involution $\varphi : C_{\eta} 
\leftrightarrow 
C_{\eta}$ over $C$ acts by $-1$.
One quickly finds that these are all in the component 
$(g-1) \otimes \c^{3}_{\rho}$,
where
$$
\c^{3}_{\rho} = \{\boldx \in \c^4 \ |\ x_1 + x_2 + x_3 + x_4 =0\}
$$
with $S_4$ permuting coordinates. The two varieties are covered, 
respectively,
by the subspaces
$$
\begin{array}{rclrcll}
H^1(R_{\eta}, \oo) &=& \c^{g-1}\otimes \boldl, &\textstyle
                          \boldl &=& {1\over 2}(3,-1,-1,-1)&\in 
\c^{3}_{\rho},\\
H^1_-(C_{\eta},\oo) &=& \c^{g-1}\otimes \boldn, & 
                          \boldn &=& (1,1,-1,-1)&\in \c^{3}_{\rho}.\\
\end{array}
$$
One notes at once (see \ref{S4notation}) that $(1+\varphi^2)\boldl = 
\boldn$
and
$(1+\tau + \tau^2)\boldn = 2 \boldl$.
This demonstrates the proposition as far as an {\it isogeny} 
$\psi: J_{R_{\eta}}\rightarrow P_{\eta}$; we have to see that $\psi$
is the Recillas isomorphism as claimed. 
But the latter, by construction, is induced from the 
inclusion $C_{\eta} \subset S^2 R_{\eta}$, and this fits into a 
commutative diagram:
$$
\begin{array}{ccc}
\gal(R_{\eta}) & \map{1+\varphi^2} & S^2 \gal(R_{\eta})\\
&&\\
\downarrow && \downarrow \\
&&\\
C_{\eta} & \hookrightarrow & S^2 R_{\eta}\\
\end{array}
$$
We conclude that via pull-back, $1+\varphi^2$ gives the polarised 
isomorphism $\psi: J_{R_{\eta}}\idfy P_{\eta}$.
\qqed

\subsection{Trigonal Schottky}

Given a choice of $\zeta \in S(\eta)/\eta$ we can identify $P_{\eta}$
with one component of $\nm^{-1}(\eta)$. When $C$ is trigonal $P_{\eta}$
is canonically identified with the Jacobian $J_{R_{\eta}}$, and the latter 
(more precisely, its Kummer) can therefore be viewed as a component
of the Schottky configuration:
$$
\sch_{\zeta}:
J_{R_{\eta}} \rightarrow \nm^{-1}(\eta) \rightarrow \su(2).
$$
In Proposition \ref{rec_bundles} 
we shall describe explicitly how rank 2 vector bundles on $C$
are constructed from line bundles on $R_{\eta}$ using diagrams 
(\ref{S4galois}) and (\ref{S4galoisdual}); 
in the next section we will
then extend the construction to rank 2 vector bundles on $R_{\eta}$.

Let us denote by $\pi_R$ the map $\gal(R_{\eta})\rightarrow R_{\eta}$; then 
for any line bundle $L\in J_{R_{\eta}}$ we consider the collection of
bundles on $\gal(R_{\eta})$ obtained by applying the 
$S_4$-action
to $L_0 := \pi_R^* L$. Since $L_0$ is $S_3$-invariant there are just 
four 
isomorphism classes, and these are in 1:1 correspondence with the left 
cosets
of $S_3 < S_4$. If we label the bundles as
$$
\begin{array}{rcl}
L_1 &=& \varphi^2 L_0 \\
L_2 &=& \varphi L_0 \\
L_3 &=& \varphi^3 L_0 \\
\end{array}
$$
then the Galois action of $S_3 < S_4$ will correspond to the usual 
permutation action on the 
subscripts $1,2,3$ (since $\varphi = (0213)$).

\begin{rem}\rm
Note that over a divisor $w+x+y+z \in g^{1}_{4} = \bp^{1}$ a point of 
the Galois curve represents an ordering of $w,x,y,z$ and the fibres 
of the bundles $L_{0},L_{1},L_{2},L_{3}$ over this point are 
canonically 
isomorphic, respectively, to the fibres $L_{w},L_{y},L_{x},L_{z}$ of $L$.
\end{rem}

We now consider the rank 2 vector bundle
$$
E'_{L,N} = NL_0 L_1 \oplus N^{-1} L_2 L_3
$$
where $N$ is a line bundle on $\gal(R_{\eta})$ which will be chosen in a 
moment.
First note that if $N$ is trivial then $E' = E'_{L,\oo}$ is invariant under
the dihedral group $\<\sigma , \varphi \>$, since $\sigma E' = L_0 L_1
\oplus L_3 L_2 = E'$ (acting trivially in the fibres at fixed points of 
$\sigma$), and $\varphi E' = L_2 L_3 \oplus L_1 L_0 = E'$ (with no fixed points 
of $\varphi$ to consider since $C$ is smooth). It follows that $E'_{L,\oo}$ 
descends to a bundle on~$C$.

More generally, the same is true of $E'_{L,N}$ provided $N$ is the pull-back 
of a line bundle in 
$$
\ker(1+\varphi) = \nm^{-1}(\oo_C)\cup \nm^{-1}(\eta)\subset J_{C_{\eta}},
$$ 
since then $N$ will be invariant
under $\sigma$ and anti-invariant under $\varphi$.

On the other hand, it follows from Proposition \ref{Hodgedecomp} that 
$\pi_R^* J_{R_{\eta}}
\subset \ker(1+\varphi+\varphi^2+\varphi^3)$; hence $\det E'_{L,N} = 
\oo_{\gal(R_{\eta})}$. Consequently the bundle $E_{L,N} \rightarrow C$ 
which pulls back to $E'_{L,N}$ satisfies 
$$
\det  E_{L,N} \in \{\oo_C,\eta\}.
$$
It turns out that $\det  E_{L,N} = \oo_C$ if and only if $N$ comes
from $\nm^{-1}(\eta)$; and that we recover the Schottky map 
$\sch_{\zeta}:
J_{R_{\eta}}
\rightarrow \su(2)$ if we choose $N$ to be the pull-back of $\zeta \in 
\nm^{-1}(\eta)[2] \cong S(\eta)/\eta$:

\begin{prop}
\label{rec_bundles}
For each $\zeta \in S(\eta)/\eta$ the following diagram commutes:
$$
\begin{array}{rcl}
&&P_{\eta} \\
&&\\
&{\rm Recillas}\nearrow\qquad\quad & \downarrow \sch_{\zeta} \\
&&\\
J_{R_{\eta}} & \rightarrow & \su(2) \\
&&\\
L &\mapsto & E_{L,\zeta} \\
\end{array}
$$ 
where $E_{L,\zeta}$ is defined by $\pi_{C_{\eta}}^* f^*E_{L,\zeta}=
E'_{L,\pi_{C_{\eta}}^*\circ f^*(\zeta)}$, for the maps
$
\gal(R_{\eta}) \map{\pi_{C_{\eta}}} C_{\eta} \map{f} C.
$.
\end{prop}

\prf
By Proposition \ref{recillasisom} the Recillas map $J_{R_{\eta}}
\rightarrow P_{\eta}$ is given, in the Jacobian of $\gal(R_{\eta})$, 
by $L \mapsto (1+\varphi^2)L = L_0L_1$. We shall view the latter as a line 
bundle on $C_{\eta}$; we then have to show that its image $\zeta \otimes 
f_*(L_o L_1) \in \su(2)$ pulls back to $E'_{L,\pi_{C_{\eta}}^*\circ f^*(\zeta)}$.
But this pull-back is
$$
\begin{array}{rcl}
\pi_{C_{\eta}}^* f^*(f_*(f^*\zeta \otimes L_0 L_1))&=& 
\pi_{C_{\eta}}^*((f^*\zeta)\otimes L_0 L_1 \oplus \varphi((f^*\zeta)\otimes 
L_0 L_1) ) \\
&=& \pi_{C_{\eta}}^*((f^*\zeta)\otimes L_0 L_1 \oplus (f^*\zeta^{-1})\otimes 
L_2 L_3) \\
&=& E'_{L,\pi_{C_{\eta}}^*\circ f^*(\zeta)}.
\end{array}
$$
\qqed

\subsection{Nonabelian Schottky}

Extending the preceding construction, we can
use the diagrams (\ref{S4galois}) and 
(\ref{S4galoisdual}) 
to construct Galois 
$\spin(8)$-bundles from rank~2 vector bundles on the Recillas curve.

Given a rank~2 vector bundle $E\rightarrow R_{\eta}$ we consider the  
bundles on $\gal(R_{\eta})$ obtained by applying the 
$S_4$-action
to $E_0 := \pi_{R}^* E$. Again, $S_3$-invariance implies that there are just 
four 
isomorphism classes in 1:1 correspondence with the left 
cosets
of $S_3 < S_4$. We write
$$
\begin{array}{rcl}
E_1 &=& \varphi^2 E_0 \\
E_2 &=& \varphi E_0 \\
E_3 &=& \varphi^3 E_0. \\
\end{array}
$$
The Galois action of $S_4$ on the the bundles $E_0,E_1,E_2,E_3$ then 
coincides with permutation of the indices.

\begin{lemm-def}
\label{rec_gal_bundles}
Fix an element $\zeta \in S(\eta)/\eta$, and denote its pull-back 
to $\gal(R_{\eta})$ by $\zeta'\in J_{\gal(R_{\eta})}[2]$. 
Given a semistable rank 2 vector bundle $E\rightarrow R_{\eta}$ with $\det E = 
\oo$, the rank 8 orthogonal bundles (where $E_{0}E_{1} = E_{0}\otimes 
E_{1}$ etc)
$$
\begin{array}{ccc}
            V_{E,\zeta} & = & \zeta'\otimes (E_{0}E_{1}\oplus E_{2}E_{3}),  \\
	S^+_{E,\zeta} & = &  \tau\zeta'\otimes (E_{0}E_{2}\oplus E_{3}E_{1}), \\
	S^{-}_{E,\zeta} & = & \tau^2\zeta'\otimes (E_{0}E_{3}\oplus E_{1}E_{2}),
\end{array}
$$
are semistable and
descend to a Galois $\spin(8)$-bundle $F_{E,\zeta} \in 
\nn_C$.
\end{lemm-def}

\prf
By Proposition \ref{spin-sum} these three bundles are the standard 
representations of a $\spin(8)$-bundle, and by construction they are 
permuted isomorphically 
by the 3-cycle $\tau \in S_{4}$ under the Galois action. 
The involution $\sigma$, on the other hand, lifts to an involution
of $V_F$ with 1-dimensional $-1$-eigenspace $\bigwedge^2 (E_2)_x$
at each fixed point $x\in \gal(R_{\eta})$. 
(So the dihedral invariance and descent to the curve $C$ of the 
previous section {\it fail} in this rank 2 case.) 
However, the triple is invariant under the Klein subgroup and 
by Lemma \ref{etale} the map $\gal(R_{\eta}) \rightarrow \gal(C)$ is
unramified. It follows that the bundles descend to $\gal(C)$ where they 
satisfy the conditions of Proposition \ref{galoisbundle}. 
\qqed

\begin{rem}\rm
Observe that if we define $R_{0} = C \cup \bp^{1}$ and extend any 
rank~2 bundle $E\rightarrow C$ to the trivial bundle on the component 
$\bp^{1}$, then the above construction reduces to that of 
Lemma~\ref{su2trigonal}.
\end{rem}

We shall write $\sch_{\zeta}: \ss\uu_{R_{\eta}}(2) \rightarrow \nn_C$ 
for the map determined by \ref{rec_gal_bundles}, and verify that this 
is consistent with the notation for Pryms:

\begin{theo}
\label{nonab-schottky}
For each $\eta \in J_C[2]$ and choice of $\zeta \in S(\eta)/\eta$
the following diagram commutes:
$$
\begin{array}{ccc}
P_{\eta}& \idfy J_{R_{\eta}} \rightarrow &\ss\uu_{R_{\eta}}(2)\\
&&\\
\sch_{\zeta}\downarrow \qquad&& \qquad\downarrow\sch_{\zeta}\\
&&\\
\ss\uu_C(2)&\longrightarrow& \nn_C\\
\end{array}
$$
\end{theo}

\prf
Let us fix the following notation for the 
various maps which we will need:
$$
\begin{array}{ccccc}
&&\gal(R_{\eta})&&\\
&\pi_{\gal(C)}\swarrow&&\searrow \pi_{C_{\eta}}\\
\gal(C)&&&&C_{\eta}\\
&\pi_C \searrow &&\swarrow f&\\
&&C&&\\
\end{array}
$$

Pick line bundles $L\in J_{R_{\eta}}$ and $N \in P_{\eta}$ 
corresponding 
under the Recillas isomorphism. $L$ gives rise to a Galois bundle 
$F_{L,\zeta} = (V_{L,\zeta},S_{L,\zeta}^+,S_{L,\zeta}^-)$ on $\gal(C)$
(via the right-hand side of the diagram) which pulls back to the 
triple 
of orthogonal bundles on $\gal(R_{\eta})$ (by \ref{rec_gal_bundles}):

\begin{equation}
\label{Lsplitting}
\begin{array}{rcl}
	\pi_{\gal(C)}^* V_{L,\zeta} & = & \zeta'\otimes\left\{
(L\oplus L^{-1})(\varphi^2 L\oplus \varphi^2 L^{-1})\oplus 
(\varphi^3 L\oplus \varphi^3 L^{-1})(\varphi L\oplus \varphi L^{-1})\right\}
,  \\
	\pi_{\gal(C)}^* S^+_{L,\zeta} & = & \tau\zeta'\otimes\left\{
(L\oplus L^{-1})(\varphi L\oplus \varphi L^{-1})\oplus 
(\varphi^2 L\oplus \varphi^2 L^{-1})(\varphi^3 L\oplus \varphi^3 
L^{-1})\right\}
, \\
	\pi_{\gal(C)}^* S^{-}_{L,\zeta} & = & \tau^2\zeta'\otimes\left\{
(L\oplus L^{-1})(\varphi^3 L\oplus \varphi^3 L^{-1})\oplus 
(\varphi L\oplus \varphi L^{-1})(\varphi^2 L\oplus \varphi^2 
L^{-1})\right\}.\\
\end{array}
\end{equation}
Via the left-hand side of the diagram, on the other hand,
$N$ gives
rise to a Galois bundle $F_{N,\zeta}\in \nn_C$ represented by a 
triple $(V_{N,\zeta},S_{N,\zeta}^+,S_{N,\zeta}^-)$ and we shall show that these bundles split 
under 
pull-back to $\gal(R_{\eta})$ into the same line bundle summands as 
(\ref{Lsplitting}).

Starting with $N \in P_{\eta}$ we construct 
$E_{N,\zeta} = \zeta\otimes f_* N \rightarrow C$; then 
$V_{N,\zeta}  =  
\c^2 \otimes \pi_C^*E_{N,\zeta} \oplus \tau(\pi_C^*E_{N,\zeta})\otimes 
\tau^2(\pi_C^*E_{N,\zeta})$ etc.
Now
$$
\begin{array}{rcl}
\pi_{\gal(C)}^*(\pi_C^* E_{N,\zeta}) &=& 
\pi_{C_{\eta}}^* f^* (\zeta\otimes f_* N) \\
           &=& \zeta'\otimes\pi_{C_{\eta}}^* (N\oplus \varphi(N))\\
           &=& \zeta'\otimes\pi_{C_{\eta}}^*(N\oplus N^{-1}),\\
\end{array}
$$
where we use the fact that $N$ is anti-invariant under the 
sheet-involution
$\varphi: C_{\eta} \leftrightarrow C_{\eta}$. It follows (denoting 
$\pi_{C_{\eta}}^* N$ by $N$ for simplicity, and using the fact that 
$\zeta' \otimes \tau\zeta' \otimes \tau^2\zeta'=\oo$)
that 

\begin{equation}
\label{Nsplitting}
\begin{array}{rcrcl}
	\pi_{\gal(C)}^* V_{N,\zeta} & = & \zeta'\otimes \left\{ \right.
2N \oplus 2N^{-1}&\oplus& (\tau N\oplus \tau N^{-1})
(\tau^2 N\oplus \tau^2 N^{-1})\left. \right\},  \\
	\pi_{\gal(C)}^* S^+_{N,\zeta} & = & \tau\zeta'\otimes \left\{ \right.
2\tau N \oplus 2\tau N^{-1}&\oplus& (\tau^2 N\oplus \tau^2 N^{-1})
(N\oplus N^{-1})\left. \right\}
, \\
	\pi_{\gal(C)}^* S^{-}_{N,\zeta} & = & \tau^2\zeta'\otimes \left\{ \right.
2\tau^2 N \oplus 2\tau^2 N^{-1}&\oplus& (N\oplus N^{-1})
(\tau N\oplus \tau N^{-1})\left. \right\}.
\\
\end{array}
\end{equation}

For each of the bundles $V_{L,\zeta},V_{N,\zeta}$ we expand the expression in 
(\ref{Lsplitting})
and (\ref{Nsplitting}). Each line bundle summand is represented by a 
vector 
in $\c^{g-1}\otimes \c^3_{\rho} \subset H^1(\gal(R_{\eta}),\oo)$.
(See the proof of Proposition \ref{recillasisom}.) 
By hypothesis $L,N$ are represented respectively by $\boldx \otimes 
\boldl, 
\boldx \otimes \boldn$ for some $\boldx \in \c^{g-1}$, where 
$\boldl = {1\over 2}(3,-1,-1,-1)$ and  $\boldn = (1,1,-1,-1)$. 
We can therefore identify 
each summand, as a vector in this representation,
using the $S_4$-action. For example, $\varphi L$
is given by the vector $\boldx \otimes \varphi \boldl = \boldx 
\otimes 
(-1,-1,3,-1)$, and $L \otimes \varphi L$ by 
$\boldx \otimes (\boldl + \varphi \boldl) = \boldx \otimes 
(2,-2,2,-2)$. 
Computing in this way, we find that the summands of 
$\pi_{\gal(C)}^* V_{L,\zeta}$ are $\boldx$ tensored with the vectors:
$$
\begin{array}{l}
\pm (1,1,-1,-1),\\
\pm (1,1,-1,-1),\\
\pm (2,-2,0,0),\\
\pm (0,0,2,-2).\\
\end{array}
$$
We now play the same game with $\pi_{\gal(C)}^* V_{N,\zeta}$, representing 
$N$ by
$\boldx\otimes (1,1,-1,-1)$, $\tau N$ by $\boldx\otimes (1,-1,1,-1)$ 
etc.
After computing we find that $\pi_{\gal(C)}^* V_{N,\zeta}$ is represented by 
{\it exactly 
the same eight vectors} as above. We conclude that $\pi_{\gal(C)}^* 
V_{L,\zeta} \cong 
\pi_{\gal(C)}^* V_{N,\zeta}$, and hence by construction  $V_{L,\zeta} 
\cong V_{N,\zeta}$. 
Repeating this calculation for the spinor bundles shows likewise that 
$S_{L,\zeta}^{\pm} \cong S_{N,\zeta}^{\pm}$, and hence $F_{L,\zeta} = 
F_{N,\zeta}$.
\qqed

\begin{rem}\rm
There is a natural `Heisenberg' action of $J_{C}[2]$ on the space 
$\nn_{C}$: for $\eta \in J_{C}[2]$ we define
$$
\begin{array}{ccc}
	V_{F} & \mapsto & \pi_C^{*}\eta \otimes V_{F},  \\
	S_{F}^{+} & \mapsto & \tau(\pi_C^{*}\eta) \otimes S_{F}^{+} , \\
	S_{F}^{-} & \mapsto & \tau^{2}(\pi_C^{*}\eta) \otimes S_{F}^{-}.
\end{array}
$$
The three image bundles are cyclically permuted by $\tau$, and the 
last two exchanged by $\sigma$, and so they define a Galois 
$\spin(8)$-bundle $\eta \cdot F \in \nn_C$. Moreover, one readily 
checks, using the relation $\pi_C^{*}\eta \otimes \tau(\pi_C^{*}\eta) 
\otimes 
\tau^{2}(\pi_C^{*}\eta)= \oo$, 
that the map $ \ss\uu_{C}(2) \hookrightarrow 
\nn_{C}$ is equivariant where $J_{C}[2]$ acts by tensor product on 
rank 2 vector bundles.
We expect that the fixed-point set of this involution is precisely 
the `Schottky locus' 
$\ss\uu_{R_{\eta}}(2) \cup \ss\uu_{R'_{\eta}}(2)$.
\end{rem}


\bigskip
\noindent
{\addressit
Department of Mathematical Sciences, Science Laboratories, South 
Road, Durham DH1~3LE, U.K.
}

\noindent
{\addressit E-mail:} {\eightrm w.m.oxbury@durham.ac.uk}\qquad
{\addressit WWW:} {\eightrm 
http://fourier.dur.ac.uk:8000/\~{}dma0wmo/}

\medskip
\noindent
{\addressit
Tata Institute of Fundamental Research, Homi Bhabha Road, Bombay 
400-005, India
}

\noindent
{\addressit E-mail:} {\eightrm ramanan@math.tifr.res.in}

\end{document}